\documentclass[12pt]{article}

\usepackage[T1]{fontenc}          
\usepackage[utf8]{inputenc}

\usepackage{amsmath,amsfonts,amssymb,amsthm,mathtools}
\usepackage{mathrsfs}             
\usepackage{bm}

\usepackage[a4paper]{geometry}
\usepackage{graphicx}
\usepackage{float}
\usepackage{placeins}
\usepackage{caption}
\usepackage{marginnote}
\usepackage{comment}

\usepackage{array}
\usepackage{multirow}
\usepackage{booktabs}
\usepackage{longtable}
\usepackage{tabularx}
\usepackage{needspace}

\usepackage[shortlabels]{enumitem}

\usepackage{algorithm}
\usepackage{algpseudocode}

\usepackage{tikz}
\usetikzlibrary{positioning}
\usetikzlibrary{calc}

\usepackage[unicode]{hyperref}
\hypersetup{
	colorlinks=true,
	linkcolor=blue,
	citecolor=blue,
	urlcolor=blue
}

\numberwithin{figure}{section}
\numberwithin{equation}{section}

\makeatletter

\newtheoremstyle{myplain}
{6pt}{6pt}%
{\itshape}
{}%
{\bfseries}
{.}%
{0.5em}%
{}%

\newtheoremstyle{mydefinition}
{6pt}{6pt}%
{\normalfont}
{}%
{\bfseries}
{.}%
{0.5em}%
{}%

\newtheoremstyle{myremark}
{6pt}{6pt}%
{\normalfont}
{}%
{\bfseries}
{.}%
{0.5em}%
{}%

\renewenvironment{proof}[1][Proof]{%
	\par\pushQED{\qed}\normalfont
	\topsep6pt\relax
	\trivlist
	\item[\hskip\labelsep\bfseries #1.]%
}{%
	\popQED\endtrivlist\@endpefalse
}

\makeatother

\theoremstyle{myplain}
\newtheorem{theorem}{Theorem}[section]

\theoremstyle{mydefinition}

\newtheorem{algo}[theorem]{Algorithm}

\theoremstyle{myremark}

\theoremstyle{myplain}
\newtheorem{theo}[theorem]{Theorem}
\newtheorem{statm}[theorem]{Statement}
\newtheorem{cor}[theorem]{Corollary}

\theoremstyle{mydefinition}
\newtheorem{defi}[theorem]{Definition}

\theoremstyle{myremark}
\newtheorem{rmrk}[theorem]{Remark}

\newcommand{\df}{\operatorname{df}}       
\newcommand{\dfc}{\operatorname{def}}     
\newcommand{\dw}{\operatorname{dw}}       
\newcommand{\Ry}{\operatorname{Ry}}       
\newcommand{\capa}{\operatorname{cap}}    
\newcommand{\ch}{\operatorname{ch}}       
\newcommand{\vcount}{\operatorname{vcount}}

\begin{document}
\title{Analysis of subsystems with rooks on a chess-board representing a partial Latin square (Part 2.)\footnote{This is the second part of our consecutive articles on Latin and partial Latin squares and on bipartite graphs. Each paper holds the terminology and notation of the previous ones~\cite{[5]}.}
}
\author{Béla Jónás}
\date{Revised August 2026 \hspace{2pt} Budapest, Hungary}

\maketitle

\begin{abstract}
A partial Latin square of order $n$ can be represented by a $3$-dimensional chess-board of
size $n\times n\times n$ with at most $n^2$ non-attacking rooks. In Latin squares, a
subsystem and its remote mate together have as many rooks as their capacity, which yields a
simple \emph{capacity condition} for completion, in fact Cruse's necessary condition for
characteristic matrices.

We prove a closed-form identity for Cruse's capacity function, from which its smallest
values follow directly. The capacity of a non-degenerated remote brick couple is at least
$n-1$, with equality precisely when one edge length is $1$ and another is $n-1$; the second
smallest value is $n$, and we list the couples attaining it. Geometrically the capacity is
affine in each variable, and its minimum level set on the cube $[1,n-1]^3$ consists of the
six edges adjacent neither to $(1,1,1)$ nor to $(n-1,n-1,n-1)$.

Andersen and Hilton, and then Andersen, listed the partial Latin squares with $n$, resp.\
$n+1$ filled cells that cannot be completed. Identifying the structures that can be
overloaded, we obtain that a PLS coming from a chess-board with at most $n+1$ rooks is
completable exactly if it satisfies the capacity condition. Finally, since the two
subsystems of a remote couple are in balance within a layer, we formulate a further
necessary condition for the completion of a layer, the \emph{balance condition}.
\end{abstract}

\textbf{MSC-Class:} 05B15


\section*{Abbreviations}
\begin{tabular}{@{}ll}
LS&Latin square\\
LSC&Latin super cube\\
$d$-LSC&$d$-dimensional Latin super cube\\
PLS&partial Latin square\\
PLSC&partial Latin super cube\\
RBC&remote brick couple\\
mRBC&minimum remote brick couple\\
RAC&remote axis couple
\end{tabular}

\section{Introduction}

Take a 3-dimensional chess-board $H_n^3$ with at most $n^2$ non-attacking rooks derived from a partial Latin square P by composition. Place a dot in each empty cell of $H_n^3$ that has a Hamming distance of at least 2 from each rook of $P$. We consider the dots as candidates, only a dot can be replaced by a rook. When replacing a dot by a rook, then we have to delete all dots that have a Hamming distance 1 from the „new” rook. This process is called \emph{automatic elimination}. A cell is called \emph{eliminated} if it has neither rooks nor dots. The structure of this chess-board with rooks and dots is called \emph{Partial Latin Super Cube}, or PLSC for short. For a PLSC we always assume that the automatic elimination is done. A file of a PLSC is an \emph{eliminated file} if each cell of the file is eliminated. Each file of a completion of a PLSC contains exactly one rook, consequently, a PLSC is not completable if it has an eliminated file.  

\begin{defi}
The number of rooks in a PLSC $P$ is denoted by $|P|$. The PLSC $Q$ is an \emph{extension} of $P$ if each rook in $P$ also is a rook in $Q$, denoted by $P \subseteq Q$. $Q$ is a \emph{completion} of $P$ if $P \subseteq Q$ and $|Q| = n^2$. $P$ is called \emph{completable} if such a $Q$ exists.
\end{defi}

\begin{defi}
The RBC $(T_0,T_3)$ is \emph{stuffed}, if $(T_0,T_3)$ has as many rooks as the capacity of $(T_0,T_3)$, that is $\capa(T_0,T_3)$. 
\end{defi}

\begin{defi}
The RBC $(T_0,T_3)$ is called \emph{overloaded}, if $(T_0,T_3)$ has more rooks than $\capa(T_0,T_3)$. 
\end{defi}

The capacity function returns the number of files in an $n$-brick for degenerated RBCs, so a degenerated RBC cannot be overloaded, otherwise an $n$-brick should have more rooks than parallel files.

\begin{defi}
The \emph{closure hull} of a set of cells $X$, denoted by $\ch(X)$, is a brick $T$ that contains all cells of $X$ and $\vcount(T)$ is minimal.
\end{defi}

\begin{defi}
A brick is \emph{perfect} if it has no dots.
\end{defi}

\begin{rmrk}
No rooks can be placed into a perfect brick.
\end{rmrk}

We differentiate two special cases of extension. 
\begin{enumerate}
\item
When we place rooks only outside $\ch(P)$:
\begin{defi}
The PLSC $P$ is \emph{embedded} in $Q$ if $Q$ is an extension of $P$ and the number of rooks in $\ch(P)$ is $|P|$. 
If $P$ is perfect, then all extensions of $P$ contain $P$ in this \emph{embedded way}.
\end{defi}

\item
When we place rooks only into the $\ch(P)$:
\begin{defi}
The PLSC $Q$ is an \emph{expansion} of $P$ if $Q$ is an extension of $P$ and \mbox{$\ch(P)=\ch(Q)$}.
\end{defi}
\end{enumerate}

\begin{figure}[!htb]
\centering\includegraphics[width=0.95\textwidth]{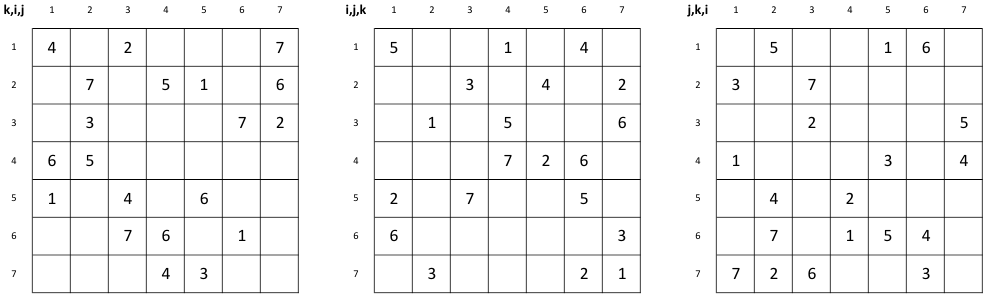}
\caption{}\label{fig1_1}
\end{figure}

If $P$ is a partial Latin square, then the conjugates of $P$ are also partial Latin squares. A PLS $P$ and its primary conjugates can be seen in the Figure~\ref{fig1_1}. For example, in the proper PLSCs, the rook and cell $(3,4,5)$ in the middle square are transformed into the cell $(5,3,4)$ in the conjugate $(k,i,j)$ (left square), and into the cell $(4,5,3)$ in the conjugate $(j,k,i)$ (right square).

Sometimes we refer to a property of the cube in Euclidean space, therefore, from now on a $3$-tuple $(a,b,c)$ has two meanings, on the one hand, it is a point in the Euclidean space, on the other hand, it identifies a cell. It will always be clear, which meaning we use. Let us consider $H_n^3$ as a cube in the Euclidean space, one vertex of the cube is at the origin and the coordinate axes are $x$, $y$ and $z$. So the vertices of $H$ in the Euclidean space are the points $(0,0,0)$, $(0,0,n)$, $(0,n,0)$, $(n,0,0)$, $(n,0,n)$, $(0,n,n)$, $(n,n,0)$ and $(n,n,n)$.

\begin{defi}
The six points $(a,b,c)$, $(c,a,b)$, $(b,c,a)$, $(a,c,b)$, $(c,b,a)$ and $(b,a,c)$ are called \emph{akin points}.
\end{defi}

If two of the coordinates are equal then we get only three distinct akin points, if all the coordinates are equal then we get only one point. Obviously, the points are equidistant from the origin in both taxicab geometry and Euclidean geometry.

\FloatBarrier

\section{Capacity Condition and its Consequences}

\begin{defi}
An arbitrary PLSC $P$ of order $n$ satisfies the \emph{capacity condition} if 
\begin{equation}\label{(201)}
c_0 + c_3 \leq \capa(T_0,T_3)
\end{equation}
for each RBC $(T_0,T_3)$ of P where $c_0$ and $c_3$ are the number of rooks in $T_0$ and $T_3$, respectively.
\end{defi}
Due to Distribution Theorem~\cite{[5]} there is no completion of $P$ if P does not satisfy the capacity condition.
The same result is described by Cruse~\cite[Theorem~2]{[3]} based on the characteristic matrix of the Latin squares.

The capacity condition is the case $k=d=3$ of a general family of necessary
conditions, which we now describe. Recall from~\cite{[5]} that a brick $T_0$
induced by $\langle E_1,\ldots ,E_d\rangle$ generates a partition of $H_n^d$ into
the $2^d$ bricks $T_I$, $I\subseteq\{1,2,\ldots ,d\}$, obtained by replacing $E_i$
by its complement for $i\in I$, and that $d(T_0,T_I)=|I|$.

\begin{statm}[Odd condition]\label{statm:odd}
Let $P$ be a completable PLSC of order $n$ in $H_n^d$, let $T_0$ be a brick and let
$I\subseteq\{1,2,\ldots ,d\}$ with $|I|$ odd. Then
\begin{equation}\label{(202)}
c_0+c_I\leq\frac{V(T_0)+V(T_I)}{n}.
\end{equation}
\end{statm}

\begin{proof}
Let $L$ be a completion of $P$ and let $c_0^*$, $c_I^*$ denote the numbers of rooks
of $L$ in $T_0$ and $T_I$. Since $|I|$ is odd, the Main Theorem of~\cite{[5]} gives
$\df(T_I)=-\df(T_0)$ in $L$, that is,
\[
c_0^*+c_I^*=\frac{V(T_0)+V(T_I)}{n}.
\]
A completion only adds rooks, so $c_0\leq c_0^*$ and $c_I\leq c_I^*$, whence
\eqref{(202)}.%
\end{proof}

For $d=3$ and $I=\{1,2,3\}$, Statement~\ref{statm:odd} is exactly the capacity
condition~\eqref{(201)}. If $|I|$ is even, then the Main Theorem yields an
\emph{equality} between the two unknown numbers $c_0^*$ and $c_I^*$ rather than an
upper bound on their sum, so it does not directly forbid anything; we do not use
such $I$ below.

The next statement shows where the genuinely new information sits.

\begin{statm}[Reduction]\label{statm:reduction}
Let $I\subseteq\{1,2,\ldots ,d\}$ with $|I|=k<d$ and let $m\notin I$. Then
condition~\eqref{(202)} for the pair $(T_0,T_I)$ in $H_n^d$ is the sum of the
corresponding $(d-1)$-dimensional conditions, taken over the layers $H_m(t)$ with
$t\in E_m$.
\end{statm}

\begin{proof}
Every layer $H_m(t)$ of a completion of $P$ is a $(d-1)$-LSC, so the restriction of
$P$ to $H_m(t)$ is a completable PLSC of dimension $d-1$. Since $m\notin I$, the
bricks $T_0$ and $T_I$ use the same set $E_m$ on the axis $t_m$. Hence, for
$t\in E_m$, their intersections with $H_m(t)$ are bricks of $H_m(t)$ at distance $k$
from each other in the partition generated by $T_0\cap H_m(t)$, while for
$t\notin E_m$ both intersections are empty. Applying \eqref{(202)} in dimension
$d-1$ to each layer $H_m(t)$ with $t\in E_m$ and adding the $|E_m|$ inequalities,
the left-hand sides add up to $c_0+c_I$ and the right-hand sides to
$\bigl(V(T_0)+V(T_I)\bigr)/n$.%
\end{proof}

\begin{cor}
For a given dimension $d$, the only condition of the family \eqref{(202)} that is
not a consequence of the conditions of dimension $d-1$ is the one with $I=\{1,2,\ldots ,d\}$,
that is, the one for a brick and its remote mate. For $d=3$ this is precisely the
capacity condition.
\end{cor}

The capacity condition also contains the condition of Ryser~\cite{[6]} as a special case.

\begin{statm}\label{statm:ryser}
Let $P$ be a PLSC of order $n$ derived from an $r\times s$ Latin rectangle, and let
$N(i)$ denote the number of occurrences of the symbol $i$ in it. If $P$ satisfies
the capacity condition, then
\[
N(i)\geq r+s-n\qquad\text{for every symbol }i.
\]
\end{statm}

\begin{proof}
Permute the symbol layers so that $i$ becomes the last one, and let $T_0$ be the
brick of size $r\times s\times (n-1)$ consisting of the first $r$ rows, the first
$s$ columns and the symbols different from $i$. Its remote mate $T_3$ has size
$(n-r)\times (n-s)\times 1$. All rooks of $P$ lie in the first $r$ rows and the
first $s$ columns, so $c_3=0$, while $c_0=rs-N(i)$. Furthermore
\[
\capa(n,r,s,n-1)=n^2-(r+s+n-1)n+\bigl(rs+s(n-1)+(n-1)r\bigr)=n+rs-r-s.
\]
The capacity condition $c_0+c_3\leq\capa(n,r,s,n-1)$ now reads
$rs-N(i)\leq n+rs-r-s$, that is, $N(i)\geq r+s-n$.%
\end{proof}

\begin{figure}[htb]
\centering
\begin{tabular}{c}
\hbox to .9\textwidth {\includegraphics[width=0.38\textwidth]{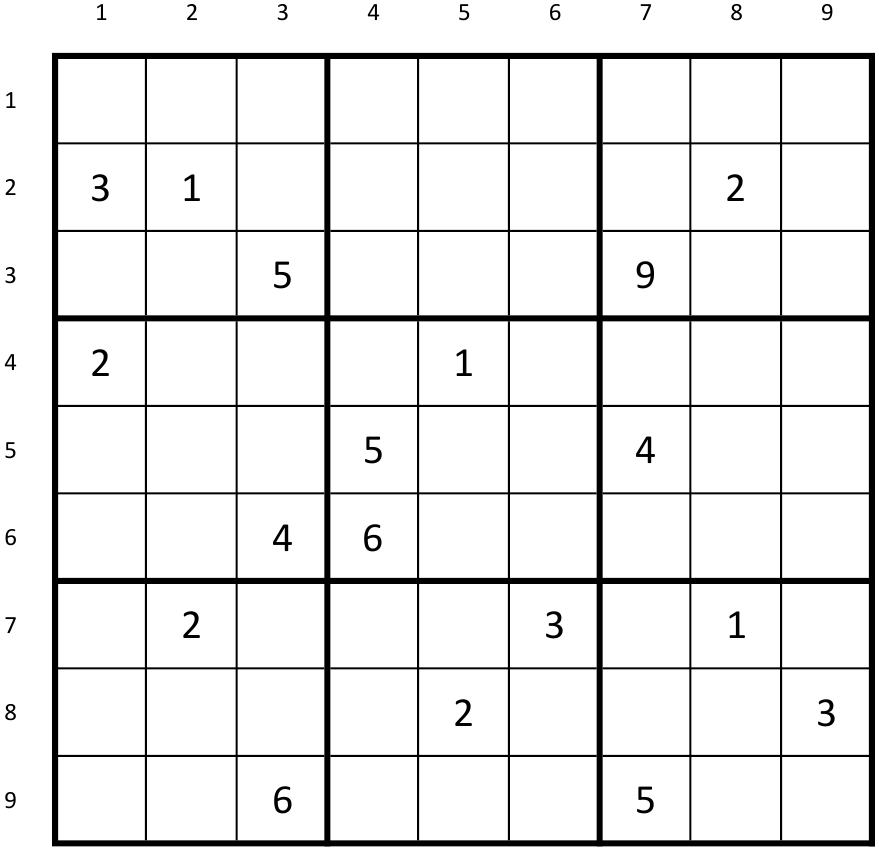}\hfill\includegraphics[width=0.38\textwidth]{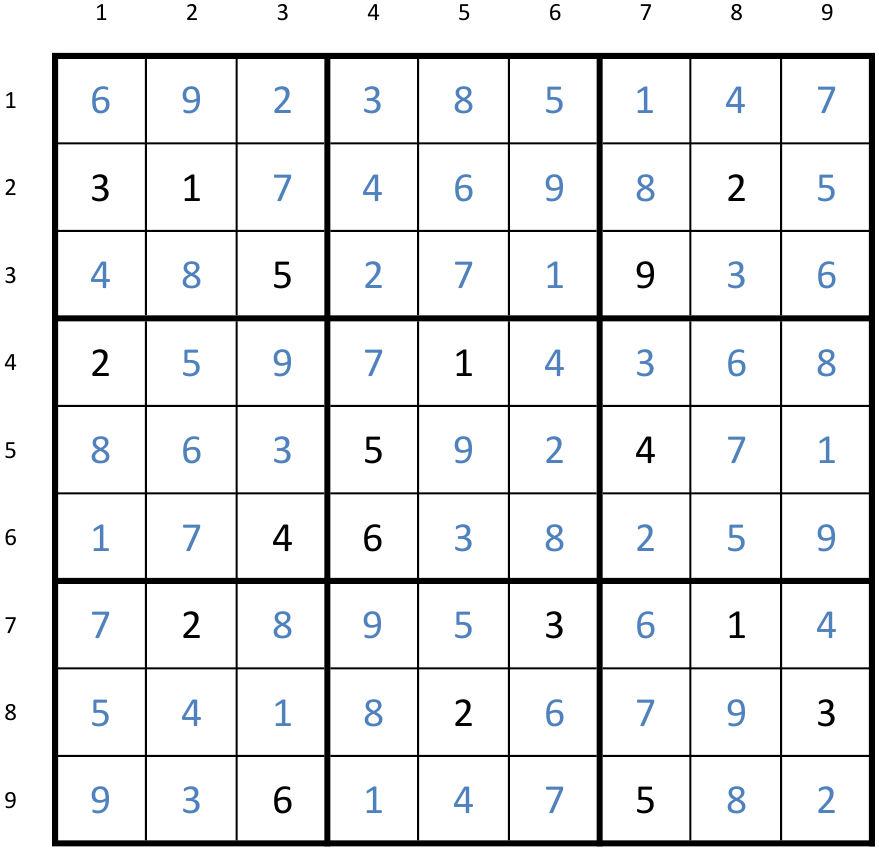}}
\end{tabular}
\caption{}\label{fig2_1}
\end{figure}

\begin{figure}[htb]
\centering
\begin{tabular}{c}
\hbox to .9\textwidth {\includegraphics[width=0.38\textwidth]{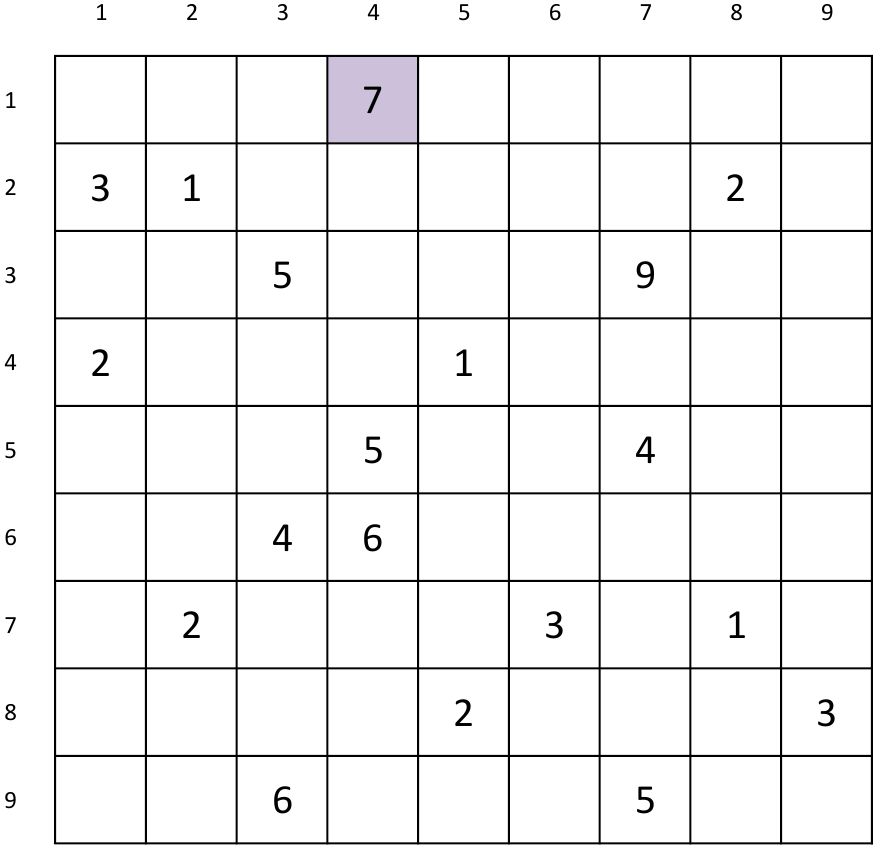}\hfill\includegraphics[width=0.38\textwidth]{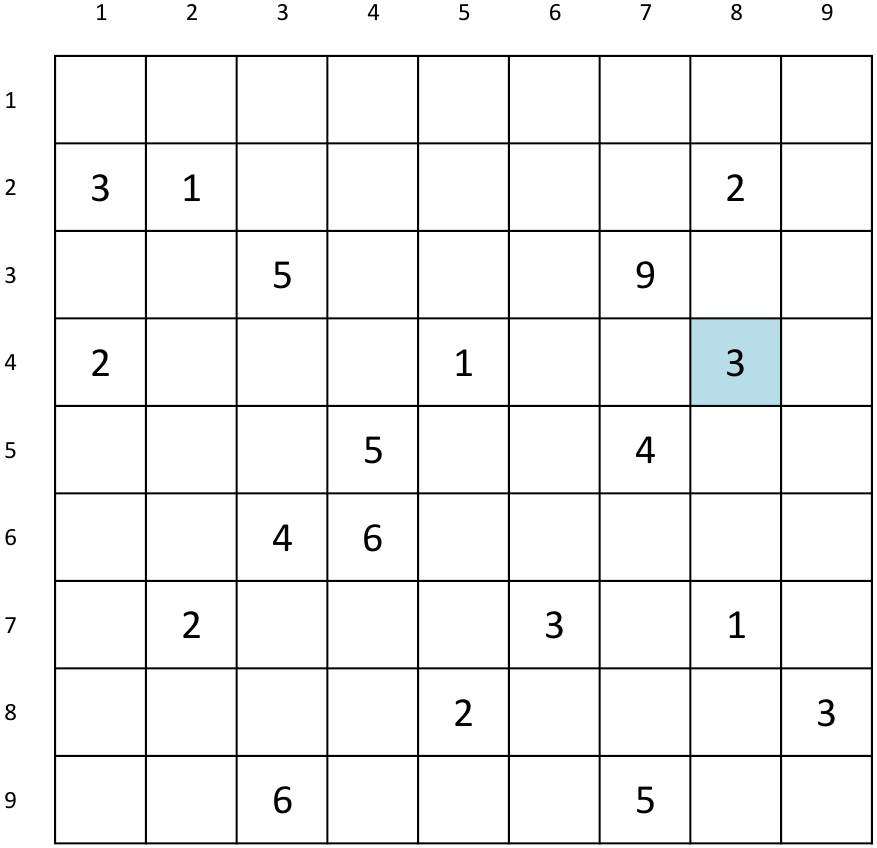}}
\end{tabular}
\caption{}\label{fig2_2}
\end{figure}

\begin{figure}[htb]
\centering
\begin{tabular}{c}
\hbox to .9\textwidth {\includegraphics[width=0.38\textwidth]{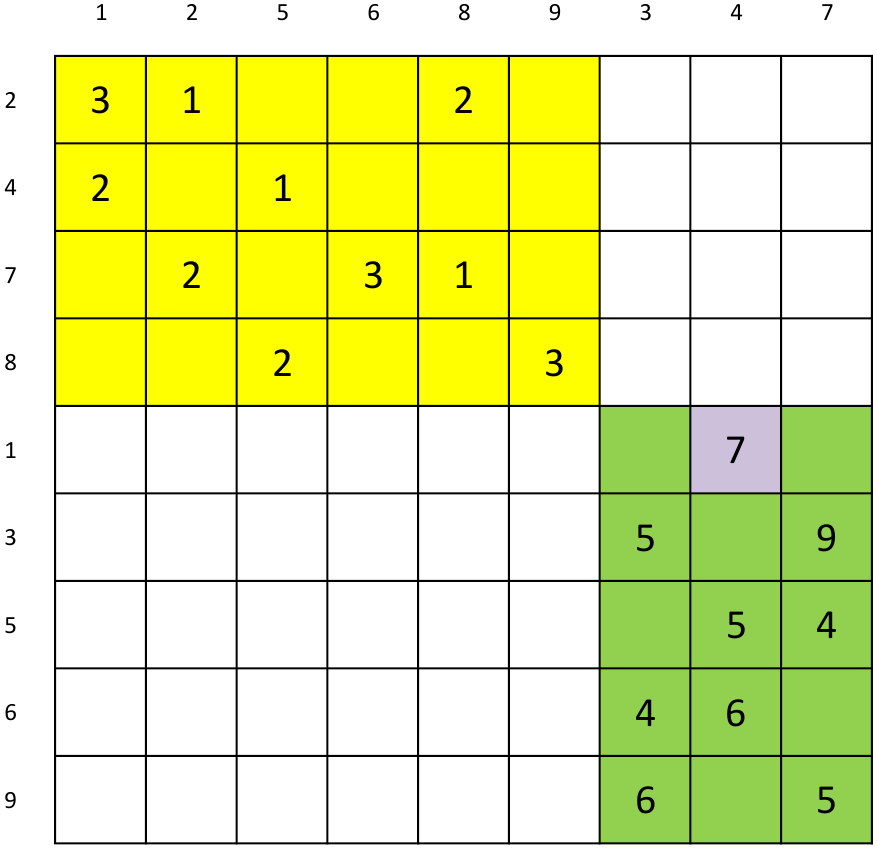}\hfill\includegraphics[width=0.38\textwidth]{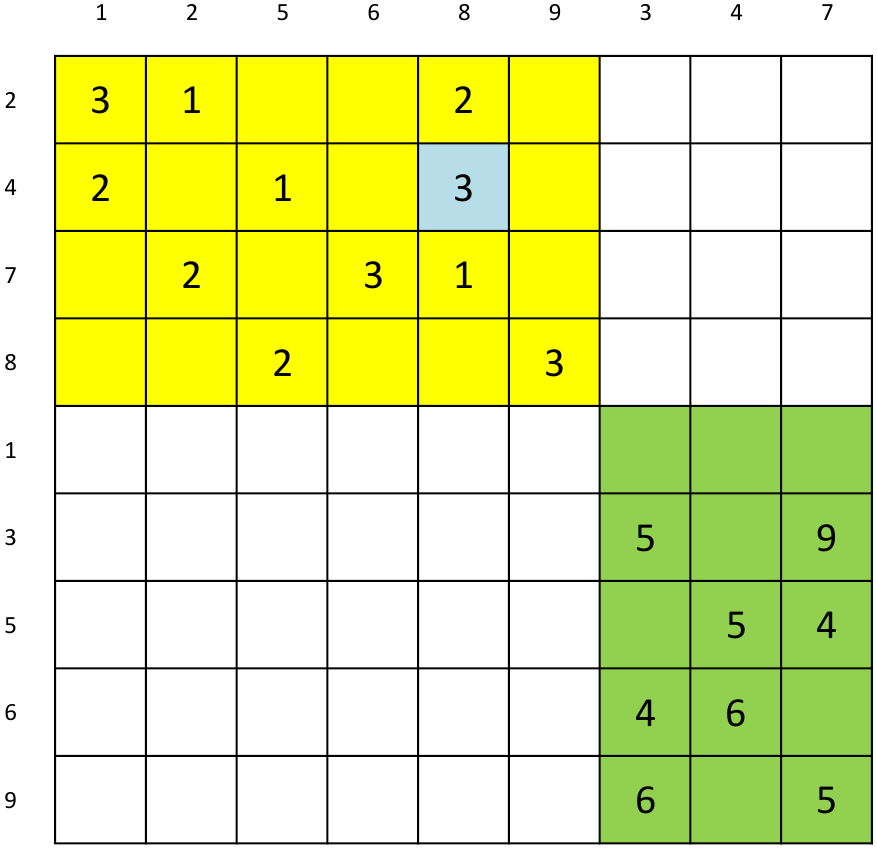}}
\end{tabular}
\caption{}\label{fig2_3}
\end{figure}

The PLSC $P$ (can also be considered as a Sudoku puzzle) in the Figure~\ref{fig2_1} has more completions. One of them is on the right-hand side of Figure~\ref{fig2_1}. If you put a 7 into the cell $(1,4)$ or a 3 into the cell $(4,8)$, the new partial Latin squares, shown in the Figure~\ref{fig2_2}, are not completable because they do not satisfy the capacity condition.
Performing the proper permutations of rows and columns we get the structures shown in the Figure~\ref{fig2_3}. The yellow brick has a height of 3, the green has a height of 6 where "height" is the number of symbol layers. The capacity of the yellow-green RBC is $(4 \times 6 \times 3 + 5 \times 3 \times 6)/9 = 18$. Some layers of the green brick have no rooks.
The RBCs have 19 rooks (symbols), so they are overloaded. Therefore, the PLSCs are not completable.

In general, rows and columns are numbered in the natural order, unless the LS or PLS is derived from another by permutation of rows and columns, in which case the original row and column coordinates allow you to restore the state before permutation. The two PLSs in the Figure~\ref{fig2_3} are derived from the PLSs in the Figure~\ref{fig2_2} by permuting the rows and columns.

\FloatBarrier

\section{Incompletable PLSCs}

\begin{figure}[htb]
\centering\includegraphics[width=0.95\textwidth]{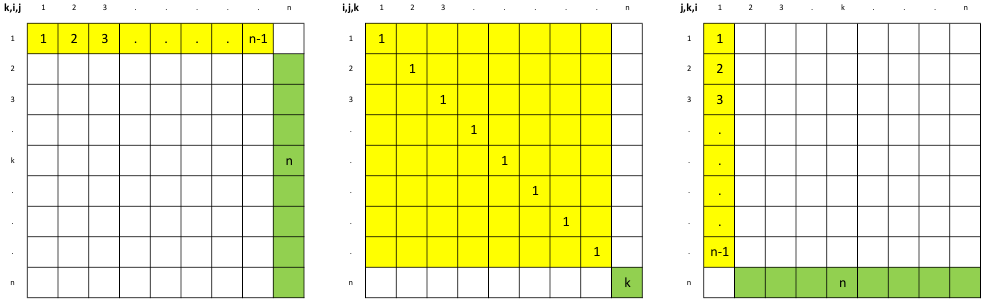}
\caption{}\label{fig3_1}
\end{figure}

\begin{figure}[htb]
\centering\includegraphics[width=0.95\textwidth]{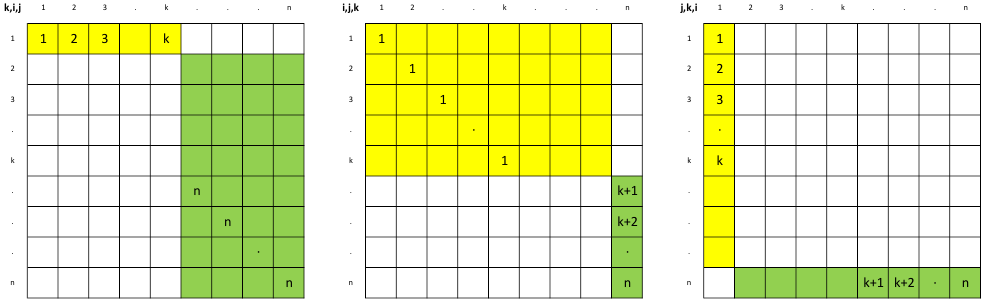}
\caption{}\label{fig3_2}
\end{figure}

Cruse~\cite{[3]} proves by induction, that the capacity condition holds, if the PLSC contains at most $n-1$  rooks. So, we need at least $n$  rooks to produce overloaded RBCs. This bound is the same as in the conjecture of Evans~\cite{[7]}, that every PLS of order $n$ with at most $n-1$ filled cells is completable, proved by Smetaniuk~\cite{[8]} and independently by Andersen and Hilton~\cite{[2]}. Here are some trivial examples in Figure~\ref{fig3_1} and Figure~\ref{fig3_2}.
The coordinates of the single green rook are ($k,n,n)$, $(n,n,k)$ and $(n,k,n)$ respectively, where $k \neq 1$. 
The capacity of the yellow-green RBCs: 
\[
\capa(n,1,n-1,n-1)=\dfrac{1 \cdot(n-1)(n-1)+(n-1)\cdot 1 \cdot 1}{n} = n-1
\]

The yellow-green RBCs have $n$  rooks, so they are overloaded.

\nopagebreak

\begin{figure}[htb]
\centering
\begin{tabular}{c}
\hbox to .95\textwidth 
{\includegraphics[width=0.3\textwidth]{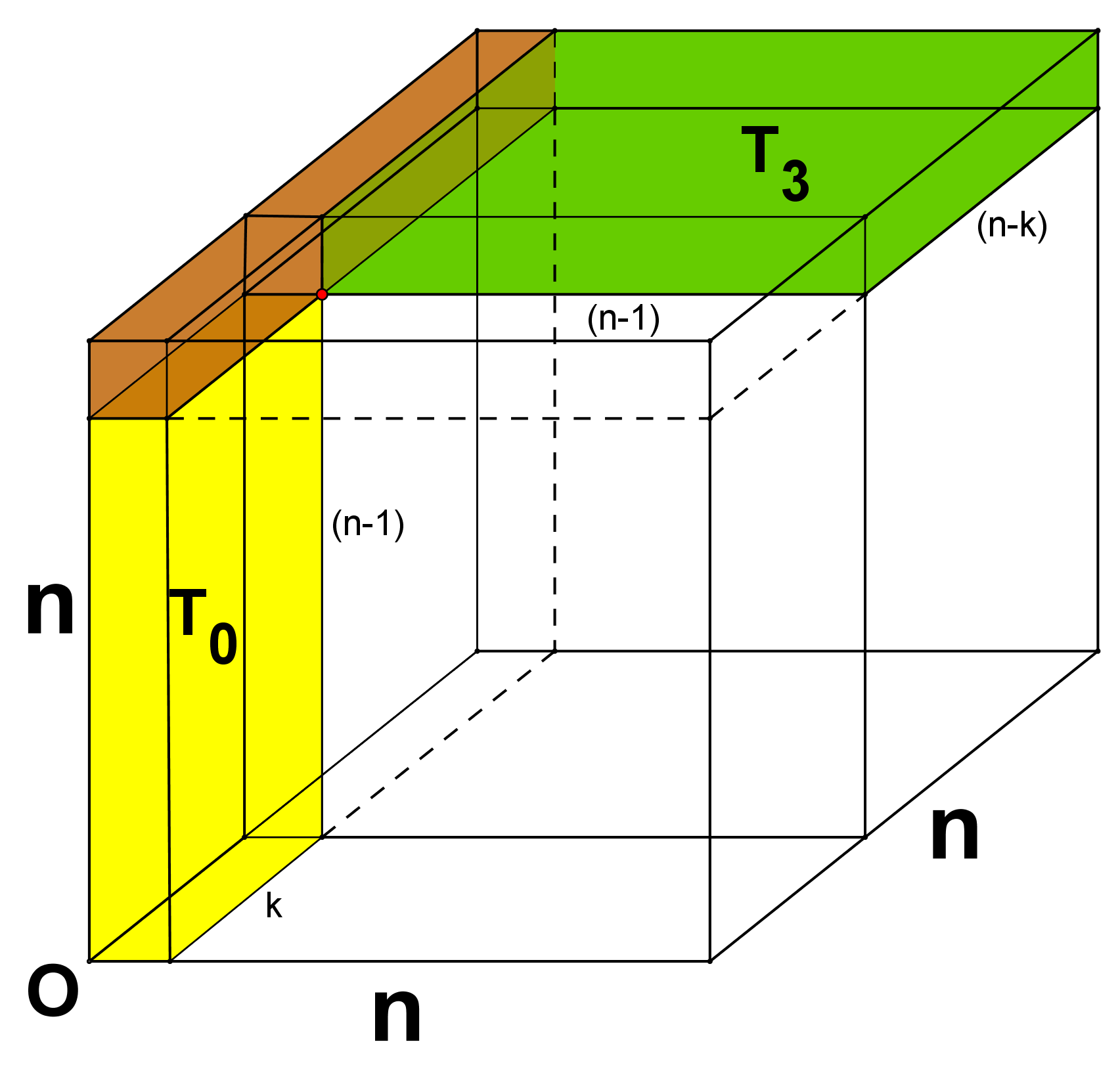}\hfill\includegraphics[width=0.3\textwidth]{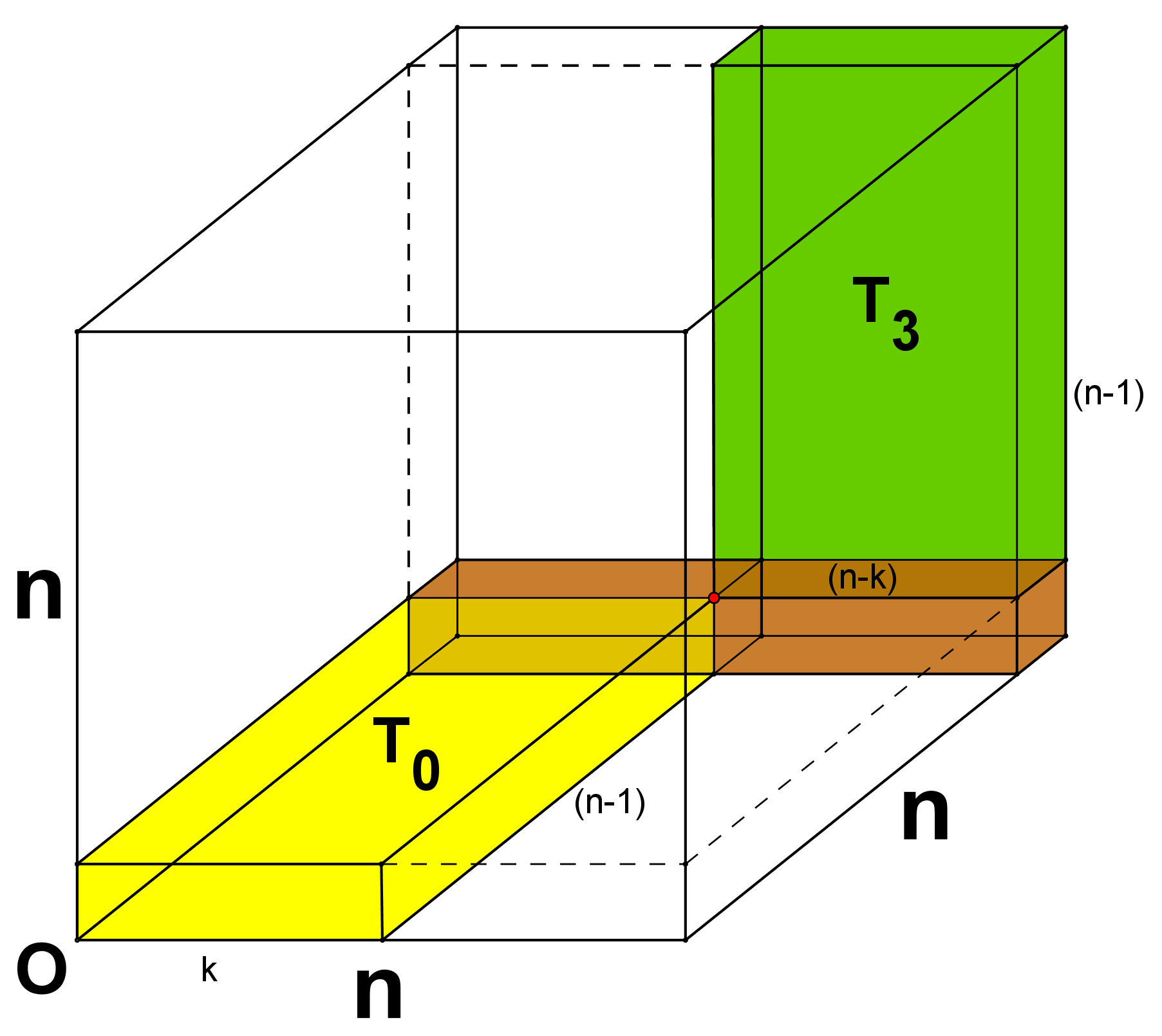}
\hfill\includegraphics[width=0.3\textwidth]{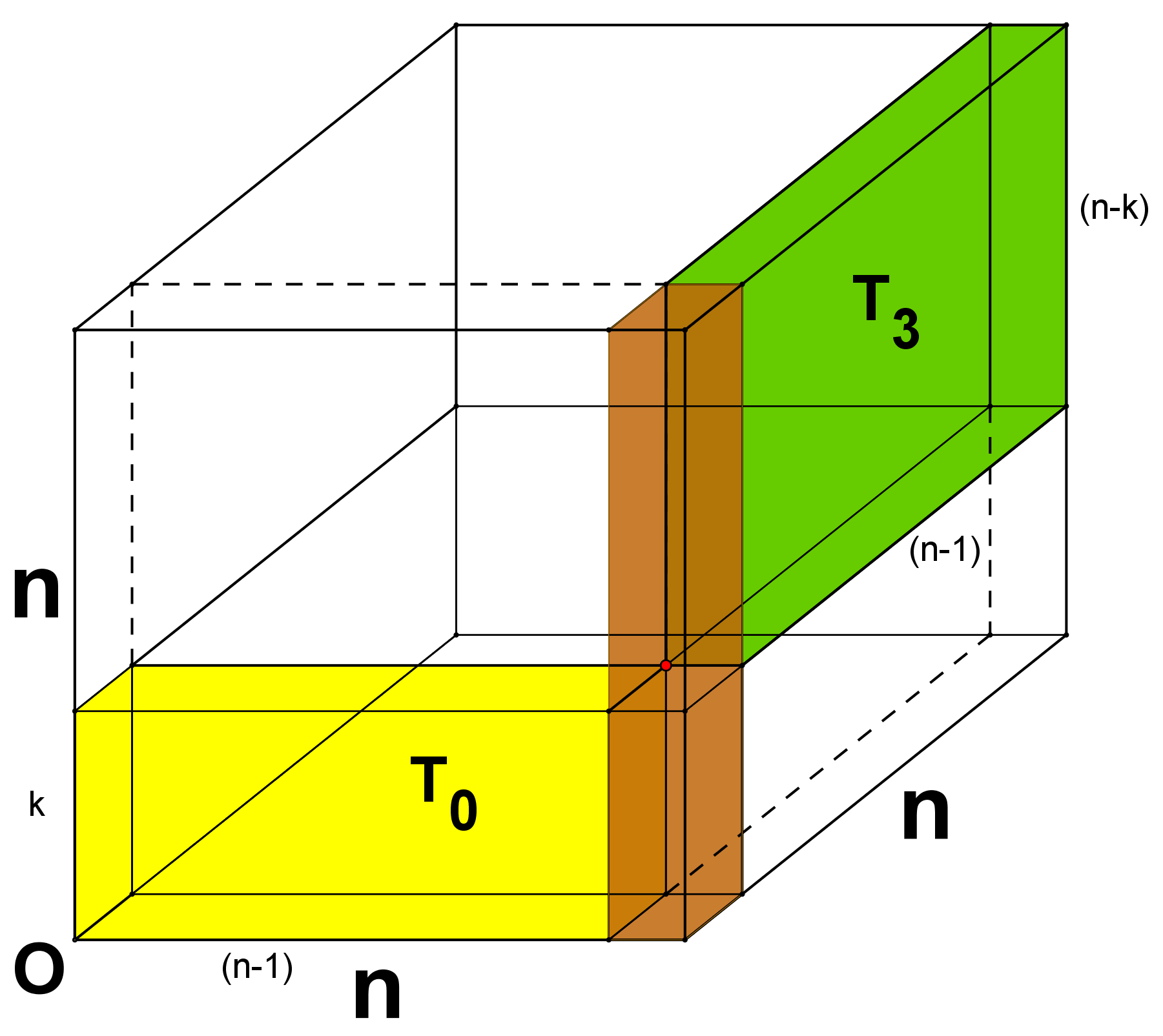}}
\end{tabular}
\caption{}\label{fig3_3}
\end{figure}

The capacity of the yellow-green RBCs:
\[
\capa(n,k,n-1,1)=\dfrac{k (n-1)\cdot 1 +(n-k)\cdot 1 \cdot (n-1)}{n} = n-1
\]
where $(1 < k < n)$.
\begin{rmrk}
The capacity of the yellow-green RBCs does not depend on $k$.
\end{rmrk}
For better visualization, first we rotated the PLSCs in the Figure~\ref{fig3_2} by $+90$ degrees and then we represented these PLSCs in a $3$-dimensional view, as shown in the Figure~\ref{fig3_3}.

The axis of the hinge is a file and the rooks of $T_0$ eliminate $k$  dots in the axis and the rooks of $T_3$ eliminate $(n-k)$ dots in the axis, so the rooks of the RBC $(T_0,T_3)$ eliminate $n$  distinct dots in this file, where $1 \leq k \leq n-1$. 
Ergo, the brown file has no dots, that means, the overloaded RBC $(T_0,T_3)$ destroys the opportunity of the completion of the brown axis.

\section{RBCs of the Smallest and Second Smallest Capacity}

Throughout this section $T_0$ is a real brick of size $a\times b\times c$, so that
$1\leq a,b,c\leq n-1$, and $T_3$ is its remote mate. We identify $T_0$ with the
integer point $(a,b,c)$ of the cube
\[
H^*=\{(a,b,c)\;:\;1\leq a,b,c\leq n-1\},
\]
whose vertices are the eight points with $a,b,c\in\{1,n-1\}$.

\begin{defi}
An RBC $(T_0,T_3)$ is called \emph{mRBC} (\emph{minimum} RBC) if there is a hinge
with leafs $T_0$ and $T_3$ whose axis consists of exactly one file. Three examples
can be seen in Figure~\ref{fig3_3}.
\end{defi}

Everything we need about the smallest values of the capacity function follows from a
single identity.

\begin{statm}[Capacity identity]\label{statm:capid}
For all integers $n\geq 2$ and $1\leq a,b,c\leq n-1$,
\begin{equation}\label{(402)}
\begin{split}
(n-2)&\bigl[\capa(n,a,b,c)-(n-1)\bigr]\\
&=(n-1-a)(n-1-b)(n-1-c)+(a-1)(b-1)(c-1).
\end{split}
\end{equation}
Both products on the right-hand side are non-negative.
\end{statm}

\begin{proof}
The capacity function is affine in each of its three variables separately, since
\[
\capa(n,a,b,c)=a(b+c-n)+\bigl[n^2-(b+c)n+bc\bigr],
\]
and $\capa$ is symmetric in $a$, $b$, $c$. Put $g(a,b,c)=\capa(n,a,b,c)-(n-1)$.
Substituting $a=1$ and $a=n-1$ gives
\[
g(1,b,c)=(n-1)^2-(n-1)(b+c)+bc=(n-1-b)(n-1-c),
\]
\[
g(n-1,b,c)=(b-1)(c-1),
\]
the second one by the central symmetry $\capa(n,a,b,c)=\capa(n,n-a,n-b,n-c)$.
For $n\geq 3$ and $1\leq a\leq n-1$ write $a=(1-\lambda)\cdot 1+\lambda\cdot (n-1)$
with $\lambda=(a-1)/(n-2)\in[0,1]$. As $g$ is affine in $a$,
\[
g(a,b,c)=(1-\lambda)\,g(1,b,c)+\lambda\, g(n-1,b,c),
\]
and multiplying by $(n-2)$ yields \eqref{(402)}. For $n=2$ both sides are $0$.
Finally, $1\leq a,b,c\leq n-1$ implies $n-1-a\geq 0$ and $a-1\geq 0$, and likewise
for $b$ and $c$, so both products are non-negative.%
\end{proof}

\begin{statm}\label{statm:mincap}
The minimum value of the capacity function for non-degenerated RBCs is $n-1$, and
$\capa(n,a,b,c)=n-1$ if and only if one of $a$, $b$, $c$ equals $1$ and another one
equals $n-1$, that is, if and only if $(a,b,c)$ defines an mRBC.
\end{statm}

\begin{figure}[htb]
\centering
\begin{tabular}{c}
\hbox to .9\textwidth {\includegraphics[width=0.44\textwidth]{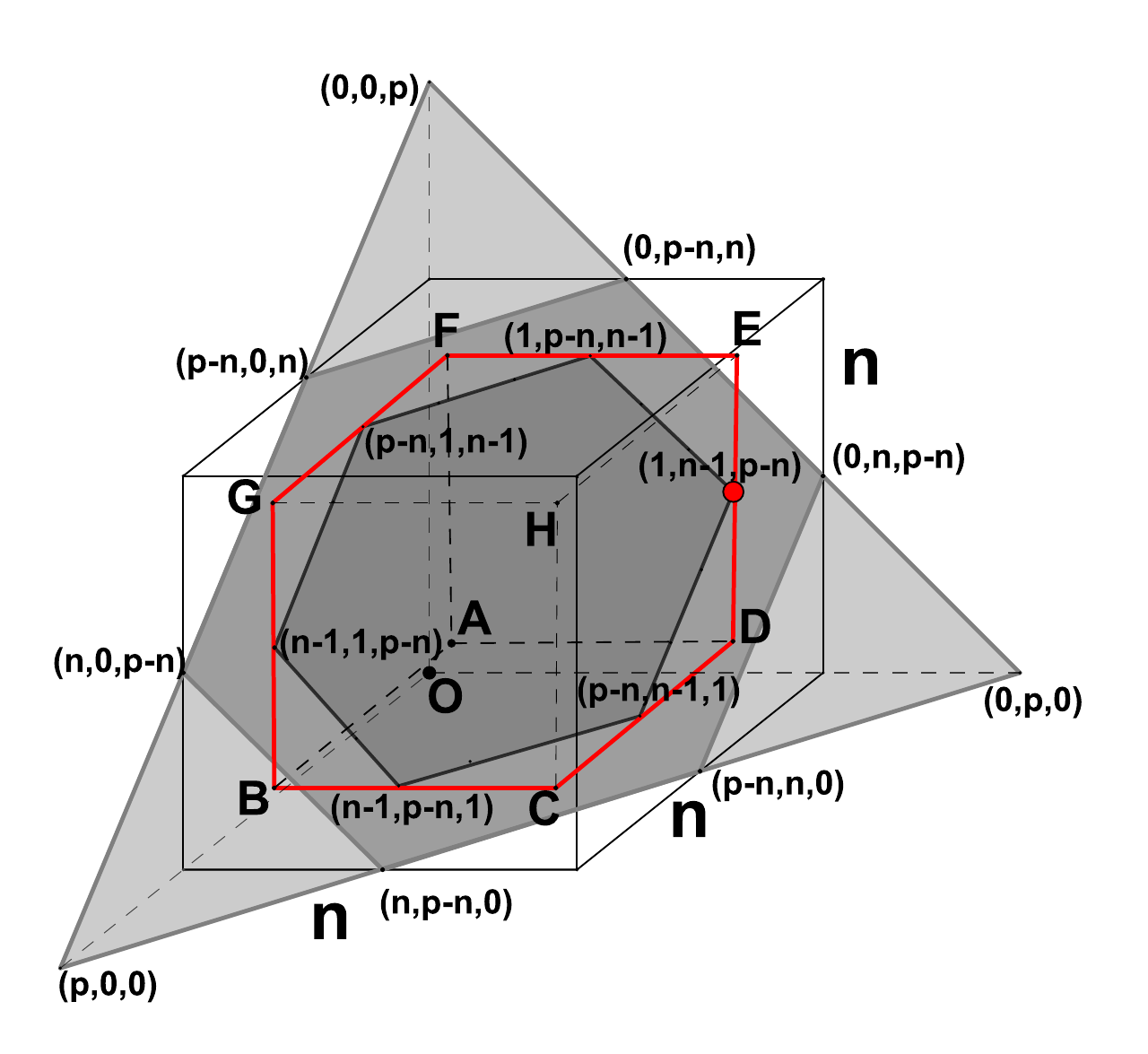}\hfill\includegraphics[width=0.44\textwidth]{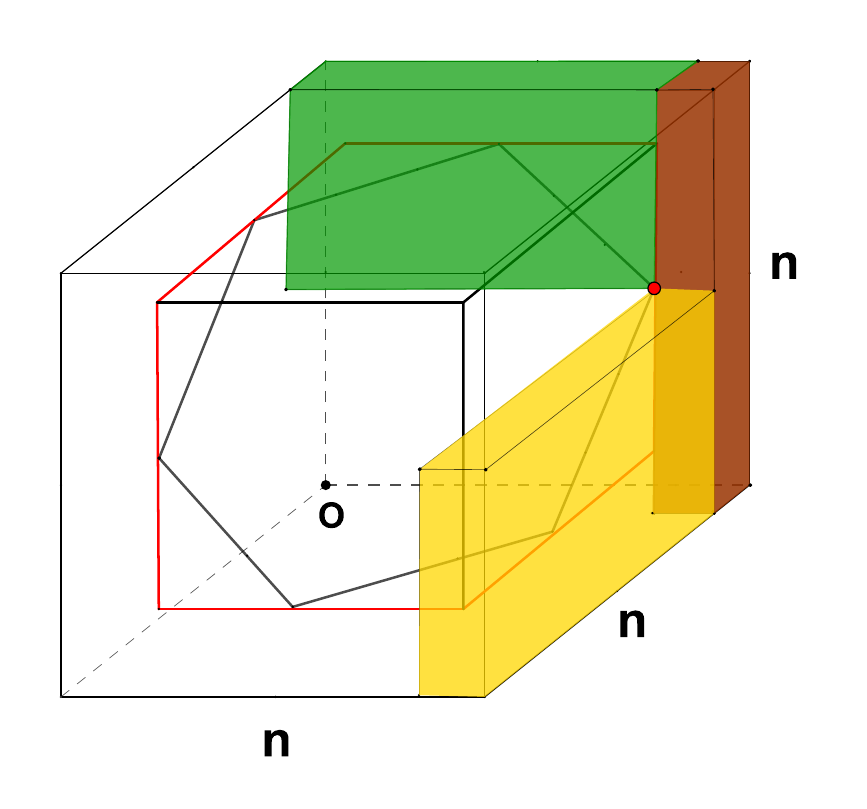}}
\end{tabular}
\caption{}\label{fig4_3}
\end{figure}

\begin{proof}
Let first $n\geq 3$, so that $n-2>0$. By Statement~\ref{statm:capid} the right-hand
side of \eqref{(402)} is non-negative, hence $\capa\geq n-1$. Equality holds exactly
when both products vanish, that is, when some coordinate equals $n-1$ and some
coordinate equals $1$. Since $1\neq n-1$ for $n\geq 3$, these cannot be the same
coordinate, so the equality set consists of the six sets
\[
\{a=1,\,b=n-1\},\quad\{a=1,\,c=n-1\},\quad\{b=1,\,a=n-1\},
\]
\[
\{b=1,\,c=n-1\},\quad\{c=1,\,a=n-1\},\quad\{c=1,\,b=n-1\},
\]
the third coordinate being arbitrary in $\{1,2,\ldots ,n-1\}$. These are precisely
the six edges of $H^*$ that are adjacent neither to the vertex $(1,1,1)$ nor to the
vertex $(n-1,n-1,n-1)$; they are the red edges on the left-hand side of
Figure~\ref{fig4_3}. A point of such an edge, say $(1,n-1,c)$, defines the hinge
whose axis is the single file lying over that edge, so it defines an mRBC, and
conversely. For $n=2$ the only admissible triple is $a=b=c=1$, and
$\capa(2,1,1,1)=1=n-1$, so the statement holds trivially.%
\end{proof}

\begin{rmrk}
The six red edges admit a pleasant geometric description. Place the outer cube
$H_n^3=[0,n]^3$ and the inner cube $H^*=[1,n-1]^3$ concentrically in Euclidean space
and cut them by the plane $x+y+z=p$. For $n+1\leq p\leq 2n-1$ this plane meets each
cube in a hexagon; the two hexagons lie in the same plane, are concentric, and have
pairwise parallel edges, the inner one being a shrunk copy of the outer. Both
degenerate to a triangle at the two ends $p=n+1$ and $p=2n-1$. The six vertices of the
inner hexagon lie on the six red edges and have integer coordinates, and as $p$ runs
from $n+1$ to $2n-1$ they sweep out the red edges. Each such vertex, being a point of a
red edge, is the leaf-corner of a hinge whose axis is a single file, that is, of an
mRBC; this is the correspondence used in the proof of Statement~\ref{statm:mincap}.
\end{rmrk}

\begin{rmrk}
Identity \eqref{(402)} also explains why the capacity is \emph{constant} along these
edges: along an edge only one variable varies, and $\capa$ is affine in each
variable, so a capacity that equals $n-1$ at both endpoints equals $n-1$ on the whole
edge.
\end{rmrk}

The second smallest value is obtained just as quickly.

\begin{statm}\label{statm:secondcap}
Let $n\geq 3$. The second smallest value of the capacity function for
non-degenerated RBCs is $n$, and $\capa(n,a,b,c)=n$ if and only if $(a,b,c)$ is a
permutation of $(n-1,2,2)$ or of $(1,n-2,n-2)$, or $n=4$ and $(a,b,c)=(2,2,2)$.
\end{statm}

\begin{proof}
Since $\capa$ is an integer, the next possible value after $n-1$ is $n$, and by
\eqref{(402)} we have $\capa=n$ if and only if
\begin{equation}\label{(403)}
X+Y=n-2,\qquad\text{where }X=\!\!\prod_{t\in\{a,b,c\}}\!\!(n-1-t),\quad
Y=\!\!\prod_{t\in\{a,b,c\}}\!\!(t-1).
\end{equation}
If $X=0$, then some coordinate equals $n-1$, say $a=n-1$; then
$Y=(n-2)(b-1)(c-1)=n-2$ forces $(b-1)(c-1)=1$, that is, $b=c=2$, giving the
permutations of $(n-1,2,2)$. The case $Y=0$ gives, by central symmetry, the
permutations of $(1,n-2,n-2)$. The case $X=Y=0$ is impossible for $n\geq 3$, since
then $X+Y=0\neq n-2$.

Suppose finally $X>0$ and $Y>0$, that is, $2\leq a,b,c\leq n-2$ (so $n\geq 4$). Put
$s=n-2$ and $x_t=t-1\in\{1,\ldots ,s-1\}$ for $t\in\{a,b,c\}$, so that
$X=\prod(s-x_t)$ and $Y=\prod x_t$. For each $t$ the product $x_t(s-x_t)$ is at
least $s-1$, its value at the endpoints of $[1,s-1]$. Hence, by the inequality
between the arithmetic and the geometric mean,
\[
X+Y\geq 2\sqrt{XY}=2\sqrt{\textstyle\prod_t x_t(s-x_t)}\geq 2(s-1)^{3/2}.
\]
Now $2(s-1)^{3/2}\geq s$ for every $s\geq 2$, with equality only for $s=2$.
Therefore \eqref{(403)} can hold only if $s=2$, that is, $n=4$, and then
$x_a=x_b=x_c=1$, that is, $a=b=c=2$; indeed $\capa(4,2,2,2)=4=n$.%
\end{proof}

\begin{figure}[htb]
\centering
\includegraphics[width=0.6\textwidth]{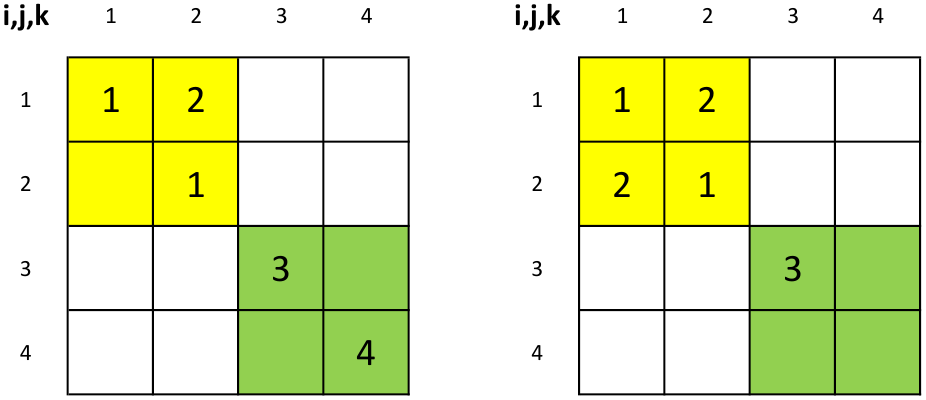}
\caption{}\label{fig4_4}
\end{figure}

For $n=3$ the two families of Statement~\ref{statm:secondcap} degenerate to the
single points $(2,2,2)$ and $(1,1,1)$. For $n=4$ they are the permutations of
$(3,2,2)$ and of $(1,2,2)$, and the extra point $(2,2,2)$ occurs as well; the
corresponding structure is the left-hand square of Figure~\ref{fig4_4}, where the
yellow and the green brick are both of size $2\times 2\times 2$ and together contain
$5=n+1$ rooks.

\begin{figure}[htb]
\centering\includegraphics[width=0.95\textwidth]{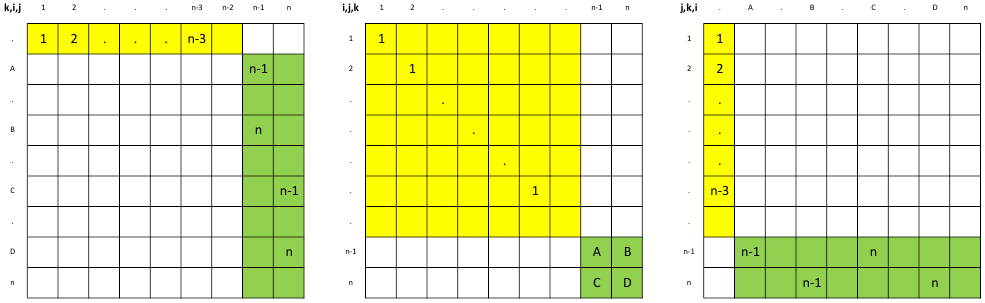}
\caption{}\label{fig4_6}
\end{figure}

The points listed in Statement~\ref{statm:secondcap} define the structures shown in
Figure~\ref{fig4_6}, or the structures derived from these by central symmetry.

We assume that the rooks derived from the symbols $A,B,C,D$ are non-attacking and are not in the first layer i.e. $A,B,C,D \neq 1$, however the cases $A = D$ and/or $B = C$ are allowed.

The capacity of the yellow-green RBCs: 
\[
\capa(n,n-2,n-2,1) = \dfrac{(n-2)(n-2)\cdot1+2\cdot2\cdot(n-1)}{n} = n
\]
The yellow-green RBCs have $(n+1)$ rooks, so they are overloaded. For $n = 4$ the square in the middle of the Figure~\ref{fig4_6} is equivalent to the second square in the Figure~\ref{fig4_4}.
There are two well-known results from Andersen-Hilton~\cite{[2]} and Andersen~\cite{[1]}. They identify and list the structures that destroy the completability of a PLS with $n$ or $n+1$ filled cells, that is, of a PLSC with $n$ or $n+1$ rooks. The listed structures are the structures shown in Figure~\ref{fig3_2} and the structures that can be obtained from the structure shown in Figure~\ref{fig4_6} by choosing the value of $A,B,C,D$ in every possible ways, and the left square in Figure~\ref{fig4_4} in case $n = 4$. We can add to this list the case $n = 3$ and the point $(1,1,1)$. Based on these results the above theorems can be restated in the following way:
\begin{theo}
Assuming $n \geq 2$, a PLSC $P$ of order $n$  with at most $(n+1)$ rooks is completable if and only if $P$ satisfies the capacity condition.
\end{theo}

\begin{figure}[htb]
\centering\includegraphics[width=0.8\textwidth]{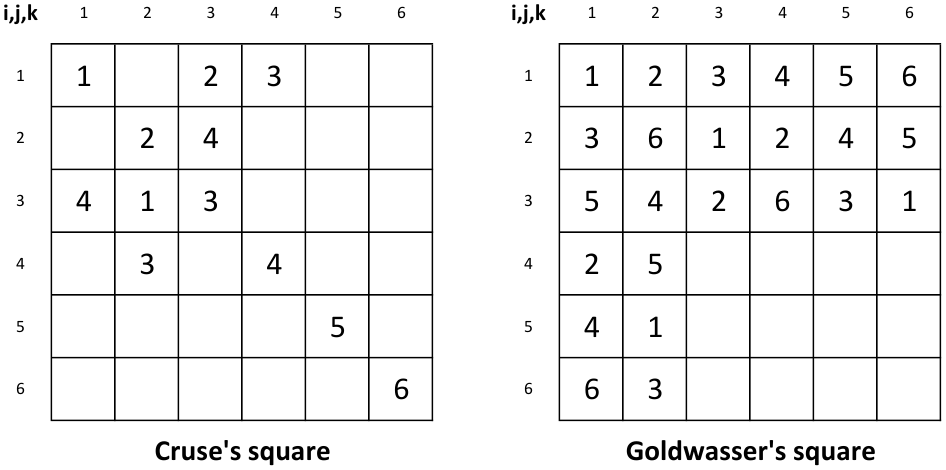}
\caption{}\label{fig4_7}
\end{figure}

Generally, the capacity condition is not sufficient for completability of a PLSC. Two iconic partial Latin squares can be seen on the Figure~\ref{fig4_7}; they are counterexamples. The first is due to Cruse~\cite{[3]}, the second, \emph{GW6}, to Goldwasser, Hilton and Patterson~\cite{[4]}, who used it in their analysis of Cropper's question.
Cruse~\cite{[3]} proved that the capacity condition holds for his square. The same method shows that \emph{GW6} satisfies the capacity condition as well, so neither of the two counterexamples is excluded by that condition.

\FloatBarrier

\section{Difference of Weights and Deficits}
\begin{figure}[htb]
\centering\includegraphics[width=0.6\textwidth]{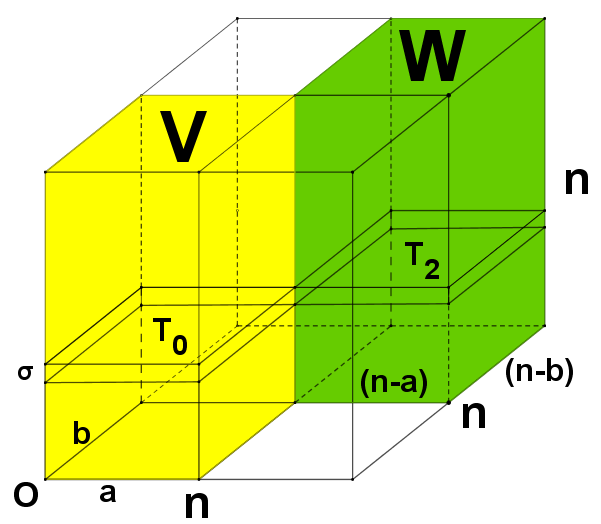}
\caption{}\label{fig5_1}
\end{figure}

\begin{defi}
The \emph{remote mate} of the axis $V$ of size $a\times b\times n$ is the axis $W$ of size $(n-a)\times (n-b)\times n$, if the edges of length $n$  of the two axes are parallel and the axes are disjoint. The pair of axes $(V,W)$ is called \emph{remote axis couple} or simply RAC if $W$ is the remote mate of $V$. 
\end{defi}
The Hamming distance between the two axes of an RAC is 2.
The intersection of $V$ and the layer $\sigma$ is denoted by $T_0$ and the intersection of $W$ and the layer $\sigma$ is denoted by $T_2$, as depicted in the Figure~\ref{fig5_1}. So, $(T_0,T_2)$ is an RBC in the layer $\sigma$ for an arbitrary $\sigma \in \{1,2,\ldots ,n\}$ and $N(\sigma,V)$ denotes the number of rooks in $T_0$ and $N(\sigma,W)$ denotes the number of rooks in $T_2$.		

\begin{defi}
The deficit of weight of the RBC $(T_0,T_2)$ is the number of rooks needed to put into $T_0$ to get the balanced state for the RBC $(T_0,T_2)$. 
So the deficit of weight of $(T_0,T_2)$ is
\[ 
\dw(T_0,T_2) = \dw(\sigma,V,W) = (a+b-n)-[N(\sigma,V)-N(\sigma,W)]
\]
In case of the balanced state $N(\sigma,V)-N(\sigma,W) = a+b-n = \Ry(T_0)$.
\end{defi}
\begin{defi}
If $(T_0,T_2)$ is an RBC of a layer then $T_0$ is \emph{underweighted} if $\dw(T_0,T_2) > 0$ and \emph{overweighted} if $\dw(T_0,T_2) < 0$.
\end{defi}
If $T_0$ is underweighted we have to put $\dw(T_0,T_2)$ rooks into $T_0$ to reach the balanced state of $(T_0,T_2)$. Since we never take rooks out of a PLSC, so we need to put $|\dw(T_0,T_2)|$ rooks into $T_2$ to reach the balanced state of RBC $(T_0,T_2)$ if $T_0$ is overweighted. 

\begin{defi}[Balance condition]\label{504}
The \emph{balance condition} holds for a layer if any perfect brick of the layer is not underweighted.
\end{defi}

We cannot put a rook in a perfect brick, so the Balance Condition \ref{504} is a necessary condition for completion of a layer.

\begin{defi}
The deficit of weight of a remote axis couple (RAC) is
\[
\dw(V,W) = \displaystyle\sum_{\sigma=1}^n (a+b-n)-[N(\sigma,V)-N(\sigma,W)] 
\]
\end{defi}
\begin{defi}
Denote $E_z(X)$ the number of files in the axis $X$ of direction $z$, that have no rooks.
\end{defi} 
Because of
\[
\dw(V,W) = n(a+b-n) -
\displaystyle\sum_{\sigma=1}^n N(\sigma,V) + 
\displaystyle\sum_{\sigma=1}^n N(\sigma,W) = 
\]
\[
[ab - \displaystyle\sum_{\sigma=1}^n N(\sigma,V)] - [(n-a)(n-b) - \displaystyle\sum_{\sigma=1}^n N(\sigma,W)]
\]
therefore,
$\dw(V,W) = E_z(V) - E_z(W)$.
Ergo, the members of a RAC are in balance, if they have the same number of empty files of direction $z$ (empty cells in the proper PLS).
If $V$ and $W$ are balanced, it does not mean that $\dw(\sigma,V,W) = 0$ for all $\sigma \in \{1,2,\ldots ,n\}$.
Obviously, $\dw(\sigma,W,V) = - \dw(\sigma,V,W)$ and $\dw(\sigma,T_2,T_0) = - \dw(\sigma,T_0,T_2)$.

Another evident condition for completion of a PLSC is the completability of all layers. If a layer of a PLSC is incompletable, then the whole PLSC is incompletable as well. Therefore, we first examine the individual layers.

\begin{figure}[htb]
\centering
\begin{tabular}{c}
\hbox to .95\textwidth {\includegraphics[width=0.44\textwidth]{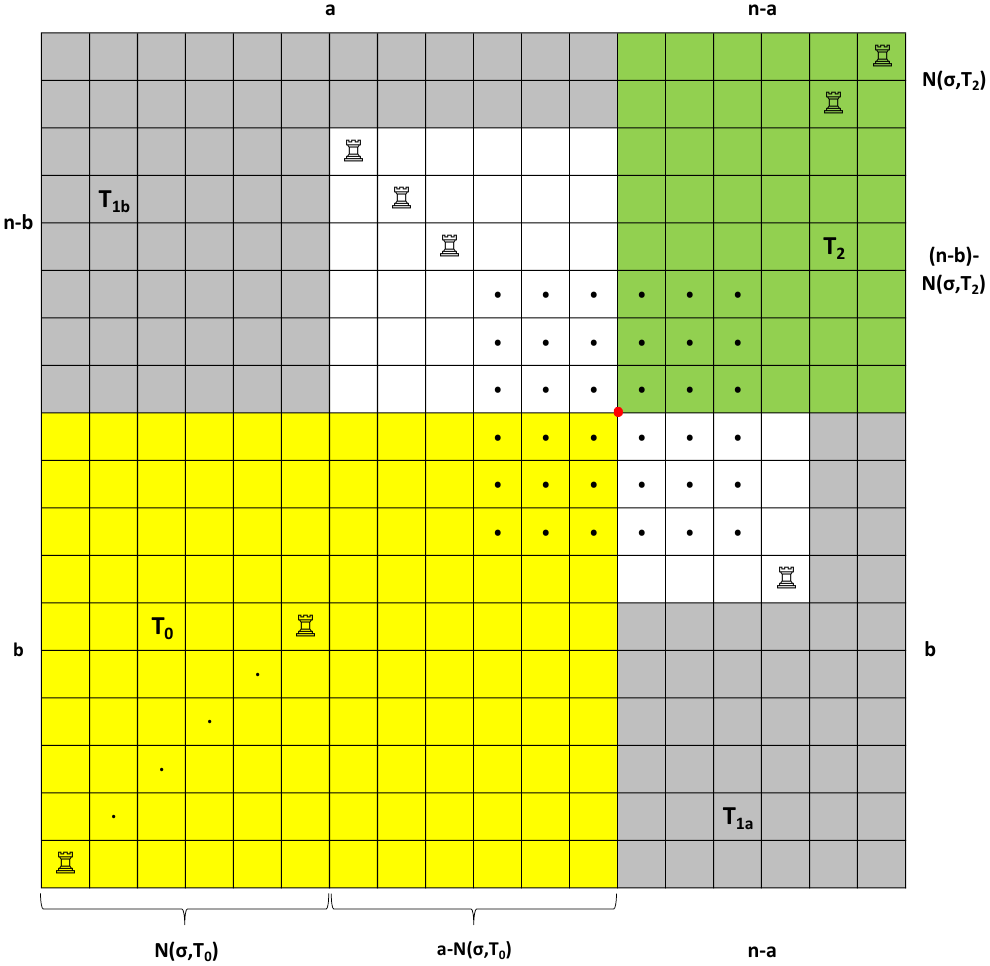}\hfill\includegraphics[width=0.44\textwidth]{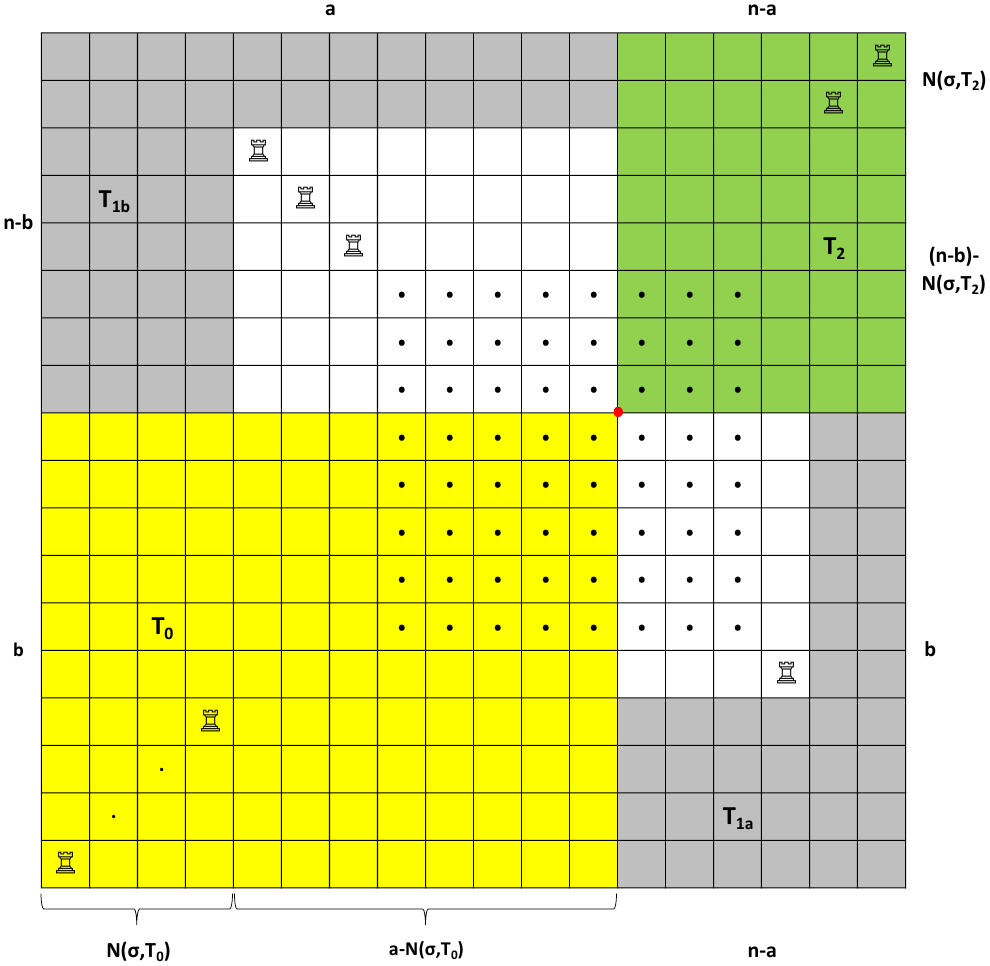}}
\end{tabular}
\caption{}\label{fig5_2}
\end{figure}

\subsection*{Case $d=2$:}	
Without loss of generality, we assume that the rows and columns have already been permuted to the form showed in the Figure~\ref{fig5_2}, so the rooks are in the corners of the brick $T_0$ and $T_2$ and in the corners of the white bricks.
Since 
\[
\dw(\sigma,T_0,T_2) = (a+b-n)-[N(\sigma,T_0)-N(\sigma,T_2)] = [a-N(\sigma,T_0]-[(n-b)-N(\sigma,T_2)]
\]
and $[a-N(\sigma,T_0)]$ is the length of the one edge of the white brick, $[(n-b)-N(\sigma,T_2)]$ is the length of the other edge of the white brick in $T_{1b}$, so the difference of the two lengths gives the deficit of weight of the RBC $(T_0,T_2)$. 
Since 
\begin{align*}
\dw(\sigma,T_0,T_2) &= (a+b-n)-[N(\sigma,T_0)-N(\sigma,T_2)]
&= [b-N(\sigma,T_0]-[(n-a)-N(\sigma,T_2)]
\end{align*}
so, we get the same result for the white brick in $T_{1a}$. This means that the RBC $(T_0,T_2)$ is balanced if and only if  the white bricks have a square shape. 
Placing a rook into a white brick leaves the deficit of weight of the RBC $(T_0,T_2)$ unchanged.
Even $\dw(T_0,T_2)$ can be calculated from the number of empty files of the white bricks.
\begin{align*}
\dw(T_0,T_2) &= (a+b-n)-(N(\sigma,T_0)-N(\sigma,T_2))\\
&= \dfrac{V(T_0)+V(T_{1b})}{n} - \dfrac{V(T_{1b})+V(T_2)}{n} - (N(\sigma,T_0) - N(\sigma,T_2))
\end{align*}

\begin{align*}
\dw(T_0,T_2) = [\dfrac{V(T_0)+V(T_{1b})}{n} - N(\sigma,T_0)] - [\dfrac{V(T_{1b})+V(T_2)}{n} - N(\sigma,T_2)]
\end{align*}

\begin{align*}
\dw(T_0,T_2) &= [\dfrac{V(T_0)+V(T_{1b})}{n} - N(\sigma,T_0) - N(\sigma,T_{1b})] \\
&- [\dfrac{V(T_{1b})+V(T_2)}{n} - N(\sigma,T_2) - N(\sigma,T_{1b})]
\end{align*}

\begin{equation}
\dw(T_0,T_2) = E(T_0\cup T_{1b}) - E(T_{1b}\cup T_2)
\end{equation} 	
In the same way we get 
\begin{equation}
\dw(T_0,T_2) = E(T_0\cup T_{1a}) - E(T_{1a}\cup T_2)
\end{equation}	
So the RBC $(T_0,T_2)$ is balanced if and only if $E(T_0\cup T_{1b}) = E(T_{1b}\cup T_2)$. In this case obviously $E(T_0\cup T_{1a}) = E(T_{1a}\cup T_2)$ as well.

\subsection*{Case $d=3$:} 
\begin{figure}[htb]
\centering\includegraphics[width=0.6\textwidth]{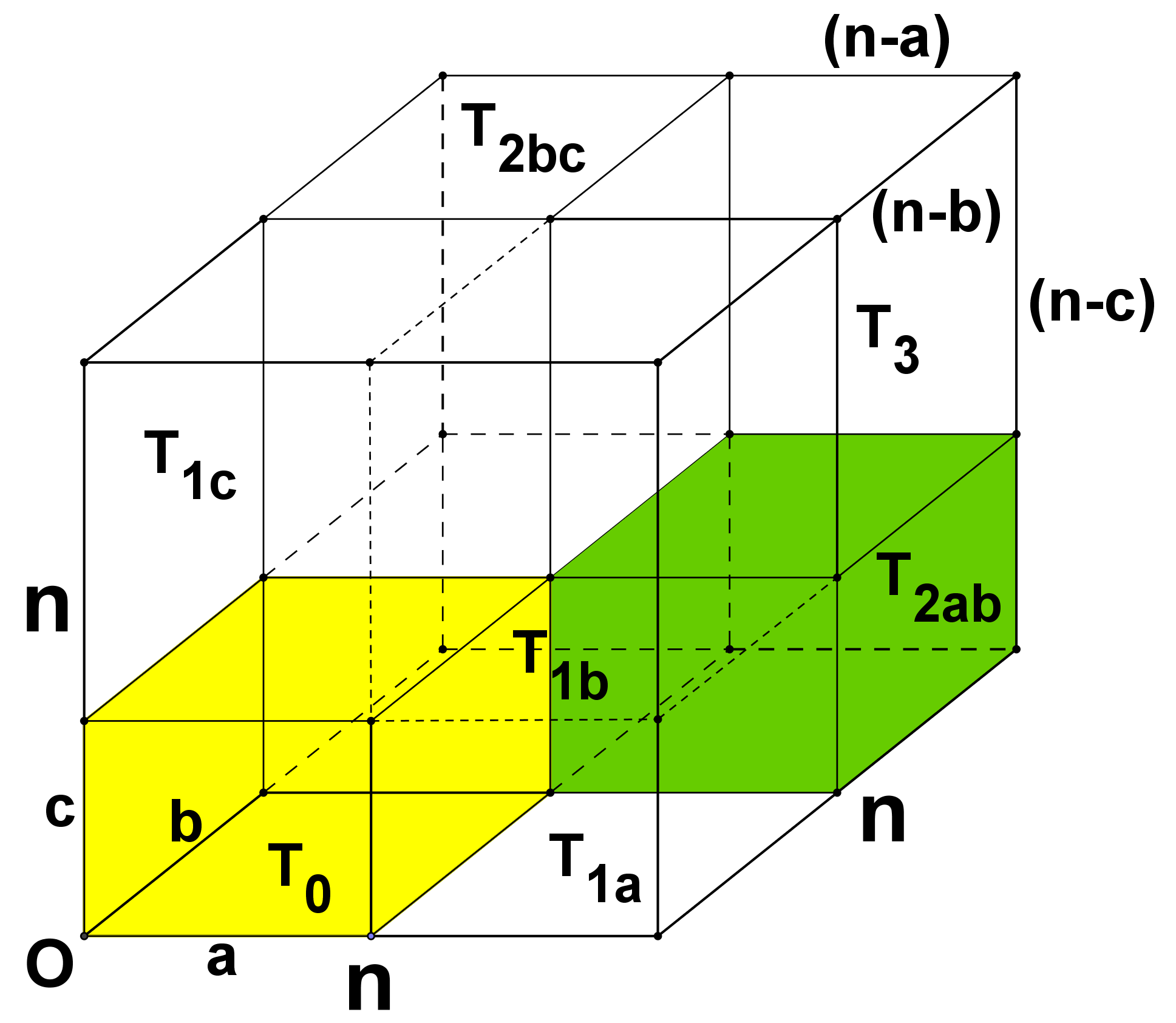}
\caption{}\label{fig5_3}
\end{figure}

Let us choose $c$ arbitrary layers in $(V,W)$. After the permutation of these layers to the bottom of the axes, these layers have the z-coordinates $1,2,\dots,c$. Let $T_0$ denote the resulting brick of size $a\times b\times c$ in $V$. The Hamming bricks generated by $T_0$ can be seen in see Figure~\ref{fig5_3}. Using the definition of the deficit of weight for a layer, we can define the deficit of weight of $(T_0,T_{2ab})$ as the sum of the deficit of weight in layer $1,2,\dots,c$ in the following way:

\begin{defi}\label{defi507}
The deficit of weight  of $(T_0,T_{2ab})$ is
\begin{equation}
\dw(T_0,T_{2ab}) = \displaystyle\sum_{\sigma=1}^c ((a+b-n)-[N(\sigma,T_0)-N(\sigma,T_{2ab})])
\end{equation}

\end{defi} 

Because
$\dw(\sigma,T_0,T_2) = E(\sigma,T_0\cup T_{1b}) - E(\sigma,T_{1b}\cup T_{2ab})$ for $\sigma  = 1,2,\dots,c$. 
Hence,
\begin{equation}\label{(501)}
\dw(T_0,T_{2ab}) = E(T_0\cup T_{1b}) - E(T_{1b}\cup T_{2ab})
\end{equation}

In the same way we get
\begin{equation}
\dw(T_0,T_{2ab}) = E(T_0\cup T_{1a}) - E(T_{1a}\cup T_{2ab})
\end{equation}

$E(T_0\cup T_{1b})$ is the number of empty files in $T_0\cup T_{1b}$ (parallel to edge $b$) and $E(T_{1b}\cup T_{2ab})$ is the number of empty files in $T_{1b}\cup T_{2ab}$ (parallel to edge $a$) and respectively $E(T_0\cup T_{1a})$ is the number of empty files in $T_0\cup T_{1a}$ (parallel to edge $a$) and $E(T_{1a}\cup T_{2ab})$ is the number of empty files in $T_{1a}\cup T_{2ab}$ (parallel to edge $b$). Therefore, $T_0$ and $T_{2ab}$ are balanced if and only if $E(T_0\cup T_{1b}) = E(T_{1b}\cup T_{2ab})$. In this case $E(T_0\cup T_{1a}) = E(T_{1a}\cup T_{2ab})$ as well.
If you stand in the middle of the white brick $T_{1a}$ and ignore the rooks in $T_{1a}$, you should see the same number of empty files when looking through the bricks $T_0$ and $T_{2ab}$.

Let $E_x(T), E_y(T)$ and $E_z(T)$ denote the number of empty files in $T$ of directions $x, y$ and $z$ respectively. So	
\[
\dw(T_0,T_{2ab}) = E_y(T_0) - E_x(T_{2ab}) = E_y(T_{1b}) - E_x(T_{1b})
\] and
\[
\dw(T_0,T_{2ab}) = E_x(T_0) - E_y(T_{2ab}) = E_x(T_{1a}) - E_y(T_{1a})
\]
If $c = n$, $V = T_0$, $W = T_{2ab}$, $Z = T_{1a}$ and $ZZ = T_{1b}$ then
\[\dw(V,W) = E_z(V) - E_z(W) = E_x(Z) - E_y(Z)\] and
\[\dw(V,W) = E_z(V) - E_z(W) = E_y(ZZ) - E_x(ZZ)\]
If you stand in the middle of the white axis $Z$ and ignore the rooks in $Z$, you should see the same number of empty files when looking through the axes $V$ and $W$. The same is true for the axes $ZZ$, $V$ and $W$ respectively.
\begin{defi}
Let $\dfc(T_0,T_3)$ denote the deficit of the RBC $(T_0,T_3)$, inclusive the degenerated cases.
\[
\dfc(T_0,T_3) = \capa(T_0,T_3) - (c_0+c_3)
\]
\end{defi}			
The deficit tells you, how many rooks need to to be put into the bricks $T_0$ and $T_3$ combined, to reach the stuffed state. If the deficit is negative, then the RBC is overloaded.
If $c = 0$, then 
\[\dfc(T_0,T_3) = \capa(T_0,W) - (c_{2ab}+c_3) = A_z(W) - (c_{2ab}+c_3) = E_z(W)
\]
If $c = n$, then 
\[
\dfc(T_0,T_3) = \capa(V,T_3) - (c_0+c_{1c}) = A_z(V) - (c_0+c_{1c}) = E_z(V)
\]

It follows from the Definition~\ref{defi507} that
\[\dw(T_0,T_{2ab}) = c(a+b-n) - (c_0-c_{2ab})
\] and so
\begin{equation}\label{(502)}
c_0 = c_{2ab}+c(a+b-n) - \dw(T_0,T_{2ab})
\end{equation}
On the other hand
\[
\dfc(T_0,T_3) = \capa(T_0,T_3) - (c_0+c_3)
\] ergo
\[
c_0+c_3+\dfc(T_0,T_3) = n^2-(a+b+c)n+ab+bc+ac
\]
that is
\begin{equation}\label{(503)}
c_0 = n^2 - (a+b+c)n+ab+bc+ac - c_3 - \dfc(T_0,T_3)
\end{equation}
From the two right sides of the equations \eqref{(502)} and \eqref{(503)} follows
\[
c_{2ab} + c(a+b-n) - \dw(T_0,T_{2ab}) = n^2 - (a+b+c)n+ab+bc+ac - c_3 - \dfc(T_0,T_3)
\] and
\begin{align*}
n^2 - (a+b+c)n+ab+bc+ac &= (n-a)(n-b) - c(n-a-b)\\ 
&= A(W) + c(a+b-n)\\
&= E(W) + c_3 + c_{2ab} + c(a+b-n)
\end{align*}
So, we get the next equation:
\begin{equation}\label{(504)}
\dw(T_0,T_{2ab}) + E(W) = \dfc(T_0,T_3)
\end{equation}
\begin{figure}[htb]
\centering
\hbox to \textwidth {\includegraphics[width=0.48\textwidth]{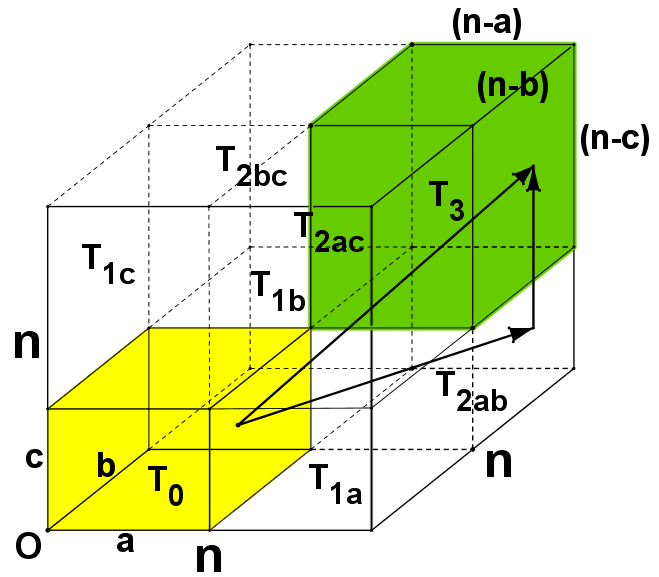}\hfill\includegraphics[width=0.42\textwidth]{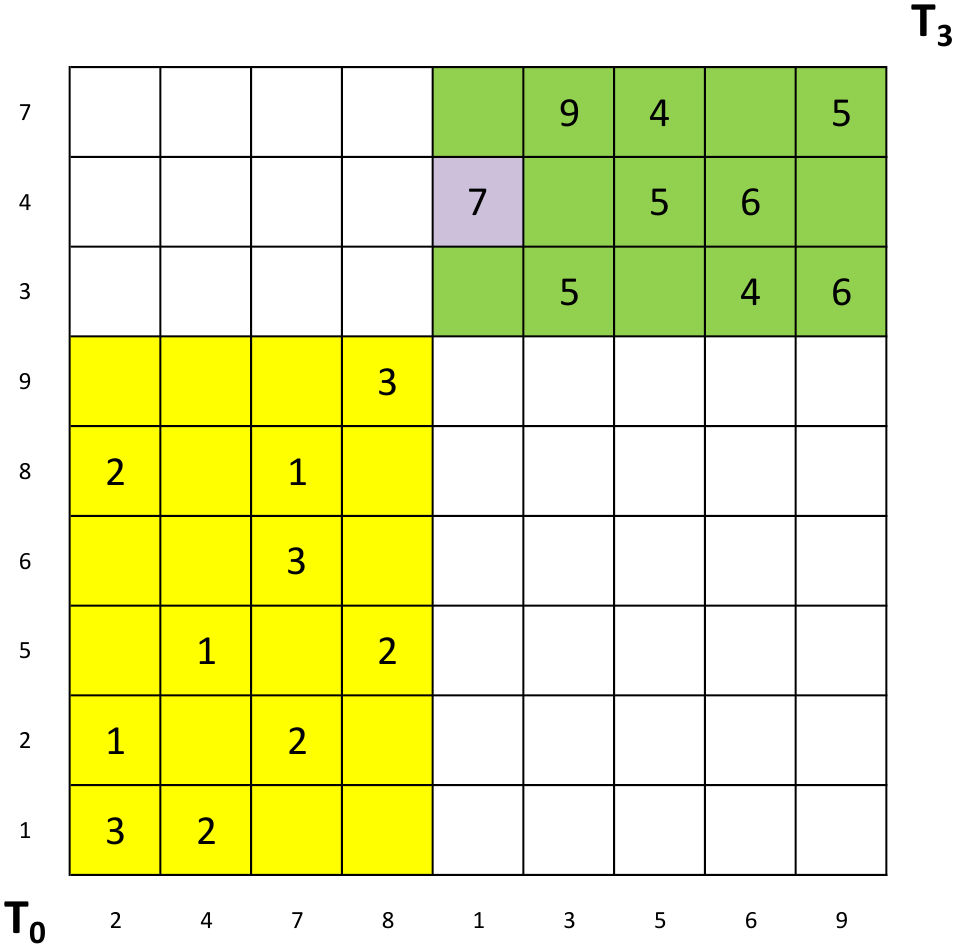}}
\caption{}\label{fig5_4}
\end{figure}

The yellow brick in the right-hand square of Figure~\ref{fig5_4} consists of $3$ layers with the following deficit of weight:
\begin{align*}
\dw(1,T_0,T_{2ab}) &= -2\\
\dw(2,T_0,T_{2ab}) &= -3\\
\dw(3,T_0,T_{2ab}) &= -2\\
\dw(T_0,T_{2ab}) &= -7\\
E(T_{2ab}\cup T_3) &= E(W) = 6\\
\dfc(T_0,T_3) &= -1
\end{align*}

We need to put $7$ rooks into the brick $T_{2ab}$ to reach the balance state. But there are only 6 empty cells in axis $W$, because we have already put the symbol 7 into the green part of a stuffed RBC. 

\begin{rmrk}
Each layer of the pair of bricks $(T_0,T_{2ab})$ can be balanced separately, you can put 2 or 3 rooks into the proper layer of brick $T_{2ab}$, but it is impossible to balance them together. This is just another explanation for the incompletability of a PLSC that has an overloaded RBC.
\end{rmrk}

\section{Capacity Condition Check and Dual Structures}
Using the equations \eqref{(501)} and \eqref{(504)} we get the next equation
\[
E(T_0\cup T_{1b}) - E(T_{1b}\cup T_{2ab}) = \dfc(T_0,T_3)- E(W) 
\]
That means
\[
\dfc(T_0,T_3) = E(W) + E(T_0\cup T_{1b}) - E(T_{1b}\cup T_{2ab}),
\]
and so
\[
\dfc(T_0,T_3) = E(T_{2ab}\cup T_3) + E(T_0\cup T_1b) - E(T_{1b}\cup T_{2ab}).
\]

The same way
\[
\dfc(T_0,T_3) = E(T_{2bc}\cup T_3) + E(T_0\cup T_{1b}) - E(T_{1b}\cup T_{2bc})
\]
If we consider the hinge with axis $T_{1b}\cup T_{2bc}$ and two leafs $T_0$ and $T_3$, we get
\begin{theo} [Hinge Deficit]
\begin{equation}\label{(505)}
\dfc(T_0\cup T_3) = E(T_{2bc}\cup T_3) + E(T_0\cup T_{1b}) - E(T_{1b}\cup T_{2bc})       	
\end{equation}
\end{theo}
\begin{proof}
\[
\begin{aligned}
&E(T_{2bc}\cup T_3) + E(T_0\cup T_{1b}) - E(T_{1b}\cup T_{2bc})\\
&= \bigl[(n-b)(n-c) - c_{2bc} - c_3\bigr] + \bigl[ac - c_0 - c_{1b}\bigr]
   - \bigl[a(n-b) - c_{1b} - c_{2bc}\bigr]
\end{aligned}
\]
and the right-hand side is equal to $n^2 -(a+b+c)n + (ab+bc+ca) - c_0 - c_3$.
\end{proof}
Due to this result, we got an algorithm to check whether a PLSC satisfies the capacity condition or not.
Let P be a PLSC of order $n$  with $|P|$ non-attacking rooks, where $|P| < n^2$.
\begin{figure}[htb]
\centering\includegraphics[width=0.6\textwidth]{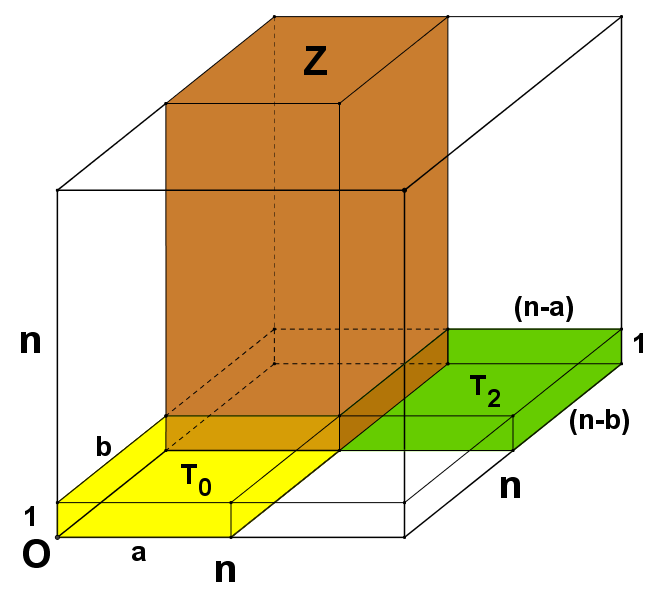}
\caption{}\label{fig5_5}
\end{figure}

Let $Z$ be a $z$-directional axis of $P$. Let $E(Z)$ be the number of $z$-directional empty files in $Z$, then $E(Z) = A(Z)-c_{1b}-c_{2bc}$. Let $E_x = 0$, $E_y = 0$. Since
\[
\dw(T_0,T_2) = (a+b-n)-\bigl[N(\sigma,T_0)-N(\sigma,T_2)\bigr]
= \bigl[a-N(\sigma,T_0)\bigr] - \bigl[(n-b)-N(\sigma,T_2)\bigr]
\]
we can count the number of empty files in the proper direction depending on the balance state of $(T_0,T_2)$, as depicted in the Figure~\ref{fig5_5}.

\begin{algo} [Capacity Condition Check]
If $\dw(T_0,T_2) < 0$, that is $T_0$ is overweighted then increment $E_y$ by $a-N(\sigma,T_0)$ otherwise increment $E_x$ by $(n-b)-N(\sigma,T_2)$ for each $\sigma \in \{1,2,\ldots ,n\}$.
\end{algo}
When all layers are done, then
\begin{align*}
E_x &= E(T_{2bc}\cup T_3) + c_{2bc}\\
E_y &= E(T_0\cup T_{1b}) + c_{1b}
\end{align*}
where the $3$-dimensional $T_0$ consists of all layers for that $E_y$ was incremented.
\[
\dfc(T_0\cup T_3) = E_x-c_{2bc} + E_y-c_{1b} - E(Z) = E_x + E_y - A(Z)
\]

It is clear from the construction that $E_x + E_y - A(Z)$ is the minimum deficit of all remote couples $(T_0,T_3)$ that form a hinge with axis $Z$.
If this value is not negative for all $z$-directional axes of $P$ then the PLSC $P$ satisfies the capacity condition.
For any remote couple $(T_0,T_3)$ there is a $z$-directional axis $Z$ so that $Z$, $T_0$ and $T_3$ form a hinge, consequently it is enough to check the number of empty files for $z$-directional axes.

Taking the advantage of the dual property of the hinge structure, we have three analogous results for volume, the number of rooks in an LSC~\cite{[5]} and the number of missing rooks in a PLSC. Following the method of this type of observations we can identify the next asymmetric correlations of bricks:
\begin{cor}[5-axes Volume]
\[
V(T_0\cup T_3) = V(T_0\cup T_{1c}) + V(T_{1b}\cup T_{2bc}) + V(T_{2ab}\cup T_3) - V(T_{1b}\cup T_{2ab}) - V(T_{1c}\cup T_{2bc})
\]
\end{cor}
and dividing both sides by $n$ we get
\begin{cor}[5-axes Capacity]
\begin{align*}
\capa(T_0,T_3) = \dfrac{V(T_0)+ V(T_{1c})}{n} &+ \dfrac{V(T_{1b}) + V(T_{2bc})}{n} +  \dfrac{V(T_{2ab})+ V(T_3)}{n}\\ &- \dfrac{V(T_{1b})+ V(T_{2ab})}{n} - \dfrac{V(T_{1c})+ V(T_{2bc})}{n}
\end{align*}
\end{cor}

\begin{theo}[5-axes Deficit]
\begin{equation}\label{(506)}
\begin{aligned}
\dfc(T_0,T_3) = E_z(T_0\cup T_{1c}) &+ E_z(T_{1b}\cup T_{2bc}) + E_z(T_{2ab}\cup T_3)\\ &- E_x(T_{1b}\cup T_{2ab}) - E_y(T_{1c}\cup T_{2bc})
\end{aligned} 
\end{equation}
\end{theo}					         
\begin{proof}
\[
\begin{aligned}
E_z(T_0\cup T_{1c}) &+ E_z(T_{1b}\cup T_{2bc}) + E_z(T_{2ab}\cup T_3)\\ &- E_x(T_{1b}\cup T_{2ab}) - E_y(T_{1c}\cup T_{2bc})\\ 
&= ab - c_0 - c_{1c} + a(n-b) - c_{1b} - c_{2bc} + (n-a)(n-b) - c_{2ab} - c_3\\
&- (n-b)c + c_{1b} + c_{2ab} - (n-a)c + c_{1c} - c_{2bc}\\ 
&= ab + a(n-b) + (n-a)(n-b) - (n-b)c - (n-a)c - c_0 - c_3
\end{aligned}
\]
and the right-hand side is equal to $n^2 -(a+b+c)n + (ab+bc+ca) - c_0 - c_3$.
\end{proof}

If $P$ is a PLSC in $H$, then the volume of $H$ is $V(H) = n^3$ and thus
$\dfrac{V(H)}{n} = n^2$. The number of missing rooks in $P$ is $n^2 - |P|$, can be denoted by $\dfc(H)$, so we can also identify the next symmetric correlations of bricks:
\begin{cor}[3-axes Volume]\label{cor606}
\[
V(T_0\cup T_3) = V(H) - V(T_{1a}\cup T_{2ab}) - V(T_{1b}\cup T_{2bc})- V(T_{1c}\cup T_{2ac})
\]
\end{cor}
\begin{cor}[3-axes Capacity]\label{cor607}
\[
\capa(T_0,T_3) = \dfrac{V(H)}{n} - \dfrac{V(T_{1a})+ V(T_{2ab})}{n} - \dfrac{V(T_{1b})+ V(T_{2bc})}{n} - \dfrac{V(T_{1c})+ V(T_{2ac})}{n}
\]
\end{cor}

\begin{theo}[3-axes Deficit]\label{theo608}
\begin{equation}\label{(618)}
\dfc(T_0,T_3) = \dfc(H) - E(T_{1a}\cup T_{2ab}) - E(T_{1b}\cup T_{2bc}) - E(T_{1c}\cup T_{2ac})
\end{equation}	
\end{theo}
\begin{proof}
\[
\begin{aligned}
\dfc(H) &= n^2 - (c_0 + c_{1a} + c_{1b} + c_2{ab} + c_{2bc} + c_{2ac} + c_3)\\
E(T_{1a}\cup T_{2ab}) &= (n-a)c + c_{1a} + c_{2ab}\\
E(T_{1b}\cup T_{2bc}) &= (n-b)a + c_{1b} + c_{2bc}\\
E(T_{1c}\cup T_{2ac}) &= (n-b)c + c_{1c} + c_{2ac}
\end{aligned}
\]
Consequently,
\[
\dfc(T_0,T_3) = n^2 - (n-a)c - (n-b)a - (n-b)c - c_0 - c_3
\]
and the right-hand side is equal to $n^2 -(a+b+c)n + (ab+bc+ca) - c_0 - c_3$.
\end{proof}

\begin{rmrk}
If $H$ has no rooks, then Theorem~\ref{theo608} is equivalent to Corollary~\ref{cor607}. If you start from the empty $H$ and replace dots by rooks in $H$ one by one, then $\dfc(H)$ decreases by one in each step, and exactly one of the four other disjoint part of the equation decreases by one as well. That is, the equation holds until you reach the position where $H$ has no dots. This is clearly not necessarily a completed LSC.
\end{rmrk}

\begin{rmrk}
The equations in \eqref{(505)}, \eqref{(506)} and \eqref{(618)} hold even if $\dfc(T_0,T_3)$ is negative.
\end{rmrk}

\FloatBarrier

\bibliographystyle{plain}

\end{document}